# LARGE DEVIATIONS FOR PROCESSES WITH DISCONTINUOUS STATISTICS

By Irina Ignatiouk-Robert

*Université de Cergy-Pontoise*

This paper is devoted to the problem of sample path large deviations for the Markov processes on $\mathbb{R}_+^N$ having a constant but different transition mechanism on each boundary set $\{x : x_i = 0$ for $i \notin \Lambda, x_i > 0$ for $i \in \Lambda\}$. The global sample path large deviation principle and an integral representation of the rate function are derived from local large deviation estimates. Our results complete the proof of Dupuis and Ellis of the sample path large deviation principle for Markov processes describing a general class of queueing networks.

**1. Introduction.** The present paper investigates sample path large deviations of Markov processes on $\mathbb{R}_+^N$ having a constant but different transition mechanism on each set

$$B_\Lambda = \{x : x_i = 0 \text{ for } i \notin \Lambda, \ x_j > 0 \text{ for } j \in \Lambda\},$$

where $\Lambda$ is a subset of $\{1, \ldots, N\}$. This property will be referred to as the *partial homogeneity of the transitions* in the following.

Such Markov processes occur in a wide class of stochastic models such as queueing networks. To establish a sample path large deviation principle in this situation, the general method of Freidlin and Wentzel [9] cannot be applied because of a discontinuity of the transition mechanism.

Our paper is motivated by various examples where a local sample path large deviation principle (see below for a precise definition) can be proved; roughly speaking, locally, the rate function can be identified by using the partial homogeneity of the processes. It is quite natural to try to extend this property in order to get a complete sample path large deviation principle.

In this paper, the problem of establishing a global principle of sample path large deviations from local large deviation estimates is investigated. It

---

is proved that, under some general conditions, such an extension holds and that the associated rate function has an integral representation.

Before formulating our results and discussing the literature of the domain, the definition of sample path large deviation principle is recalled.

For $x \in \mathbb{R}^N$, let $(X(t,x))$ be a Markov process on $E \subset \mathbb{R}^N_+$ with a given transition kernel and initial state $X(0,x) = x$. For $n \geq 1$,

$$Z_n(t,z) = X(nt, nz)/n$$

is a rescaled Markov process defined on $\mathcal{E}_n = \frac{1}{n}E$ and having initial state $Z_n(0,z) = z \in \mathcal{E}_n$. Throughout the paper, with a slight abuse of notation, we write $(X(t))$ and $(Z_n(t))$ instead of $(X(t,x))$ and $(Z_n(t,z))$. The subscript $z$ of $\mathbb{P}_z(Z_n \in \cdot)$ refers to the initial state of $(Z_n(t))$.

1.1. *Sample path large deviation principle.* The functional

$$I_{[0,T]} : D([0,T], \mathbb{R}^N_+) \to [0, +\infty]$$

is called a *rate function* if it is lower semicontinuous. The family of rescaled Markov processes $(Z_n(t), t \in [0,T])$ is said to satisfy the *sample path large deviation principle* with a rate function $I_{[0,T]}$ if:

(i) for any $x \in \mathbb{R}^N_+$ and for any open set $\mathcal{O} \subset D([0,T], \mathbb{R}^N_+)$,

$$(1) \quad \lim_{\varepsilon \to 0} \liminf_{n \to +\infty} \inf_{z \in \mathcal{E}_n : |z-x| < \varepsilon} \frac{1}{n} \log \mathbb{P}_z(Z_n \in \mathcal{O}) \geq - \inf_{\phi \in \mathcal{O} : \phi(0)=x} I_{[0,T]}(\phi);$$

(ii) for any $x \in \mathbb{R}^N_+$ and for any closed set $F \subset D([0,T], \mathbb{R}^N_+)$,

$$(2) \quad \lim_{\varepsilon \to 0} \limsup_{n \to +\infty} \sup_{z \in \mathcal{E}_n : |z-x| < \varepsilon} \frac{1}{n} \log \mathbb{P}_z(Z_n \in F) \leq - \inf_{\phi \in F : \phi(0)=x} I_{[0,T]}(\phi),$$

where $D([0,T], \mathbb{R}^N_+)$ is the set of all functions from $[0,T]$ to $\mathbb{R}^N_+$ which are right-continuous and have left limits. The set $D([0,T], \mathbb{R}^N_+)$ is endowed by the Skorohod metric.

Inequalities (1) and (2) are usually called lower and upper large deviation bounds, respectively.

A general upper large deviation bound has been obtained for the processes with discontinuous statistics by Dupuis, Ellis and Weiss in [8]. This upper bound is usually quite rough: results obtained by Alanyali and Hajek [1], Blinovskiĭ and Dobrushin [3] and Ignatiouk [11] show that the lower large deviation bound with the same rate function fails in general.

1.2. *Local sample path large deviation principle.* A *local sample path large deviation principle* with a rate function $J_{[0,T]}$ is said to hold when the following inequalities are satisfied:

$$(3) \quad \lim_{\delta \to 0} \lim_{\varepsilon \to 0} \liminf_{n \to \infty} \inf_{z \in \mathcal{E}_n : |z - \psi(0)| < \varepsilon} \frac{1}{n} \log \mathbb{P}_z(\|\psi - Z_n\|_\infty < \delta) \geq -J_{[0,T]}(\psi),$$

$$(4) \quad \lim_{\delta \to 0} \limsup_{n \to \infty} \sup_{z \in \mathcal{E}_n : |z - \psi(0)| < \delta} \frac{1}{n} \log \mathbb{P}_z(\|\psi - Z_n\|_\infty < \delta) \leq -J_{[0,T]}(\psi),$$

for *every continuous piecewise linear function* $\psi \colon [0,T] \to \mathbb{R}_+^N$. Because of the Markov property, a local sample path large deviation principle holds if the above inequalities are satisfied for any function $\psi$ having a constant velocity.

For Markov processes associated to queueing networks, a local sample path large deviation principle has been established by Dupuis and Ellis [6]. For such Markov processes, the rate function $J_{[0,T]}$ has an integral form

$$J_{[0,T]}(\psi) = \int_0^T L(\psi(t), \dot{\psi}(t)) \, dt$$

for every continuous piecewise linear function $\psi$. The *local rate function $L$* is defined by the limits

$$L(x,v) = -\lim_{T \to 0} \lim_{\delta \to 0} \lim_{\varepsilon \to 0} \liminf_{n \to \infty} \inf_{z : |z-x| < \varepsilon} \frac{1}{nT} \log \mathbb{P}_z\left(\sup_{t \in [0,T]} |x + vt - Z_n(t)| < \delta\right)$$

$$= -\lim_{T \to 0} \lim_{\delta \to 0} \limsup_{n \to \infty} \sup_{z : |z-x| < \delta} \frac{1}{nT} \log \mathbb{P}_z\left(\sup_{t \in [0,T]} |x + vt - Z_n(t)| < \delta\right)$$

and satisfies the following properties:

(a) for any $x \in \mathbb{R}_+^N$, the function $v \to L(x,v)$ is convex,

(b) for any $v \in \mathbb{R}_+^N$, the mapping $x \to L(x,v) \equiv L_\Lambda(v)$ is constant on each set $B_\Lambda = \{x : x_i = 0 \text{ for } i \notin \Lambda, x_j > 0 \text{ for } j \in \Lambda\}$.

Borovkov and Mogul′skiĭ [4] obtained a local sample path large deviation principle for partially homogeneous Markov chains with values in $\mathbb{R}_+^2$. An explicit expression for the local rate function has been derived in several situations. In [10], an explicit representation of the local rate function was obtained for Jackson networks by using the classical method of exponential change of measure and the explicit representation of the related fluid limits. Atar and Dupuis [2] give the local rate function for a class of networks for which the associated Skorohod problem has some regularity properties. Delcoigne and de La Fortelle [5] expressed the local rate function for some polling systems. In [11], the local rate function of a general class of Markov chains with discontinuous statistics was represented in terms of convergence parameters of a family of matrices.

The local large deviation principle with a rate function $J_{[0,T]}$ implies the *lower* large deviation bound (1) with the rate function $I_{[0,T]}$ defined as a lower semicontinuous regularization of the function $J_{[0,T]}$: for any $\phi \in D([0,T], \mathbb{R}^N_+)$,

$$(5) \qquad I_{[0,T]}(\phi) = \lim_{\delta \to 0} \inf_{\psi} J_{[0,T]}(\psi),$$

where the infimum is taken over all piecewise linear functions $\psi \colon [0,T] \to \mathbb{R}^N_+$ with $d_S(\phi, \psi) < \delta$ and $d_S(\cdot, \cdot)$ denotes the Skorohod metric (see Theorem 4.3 of [6]).

As it stands, the local large deviation principle is not sufficient to imply the *upper* large deviation bound (2). In this setting, a proof of the upper large deviation bound has been proposed by Dupuis and Ellis [6]. It is not entirely correct for the following reason: relation (4) shows that, for any $\varepsilon > 0$ and for any continuous piecewise linear function $\phi \colon [0,T] \to \mathbb{R}^N_+$, there exists $\delta = \delta_\phi > 0$ depending on $\varepsilon$ and also on $\phi$ such that

$$\limsup_{n \to \infty} \sup_{z \in \mathcal{E}_n \,:\, |z - \phi(0)| < \delta} \frac{1}{n} \log \mathbb{P}_z(\|\phi - Z_n\|_\infty < \delta) \leq -I_{[0,T]}(\phi) + \varepsilon.$$

Let $\mathcal{S}$ be the set of all continuous piecewise linear paths $\phi \colon [0,T] \to \mathbb{R}^N_+$. To obtain the upper large deviation bound, the arguments of [6] consist in covering a compact subset $K \subset D([0,T], \mathbb{R}^N_+)$ by a finite family of open sets $\{\psi : \|\psi - \phi\| < \delta_\phi\}$ with $\phi \in \mathcal{S}$. While the set $\mathcal{S}$ is dense in the Skorohod space $D([0,T], \mathbb{R}^N_+)$, such a covering *does not necessarily exist in general* because the quantity $\delta_\phi$ depends on $\phi$.

Moreover, (5) gives only an implicit description of the rate function $I_{[0,T]}$. Even if the closed form expression of the local rate function $L(\cdot, \cdot)$ is known, it is not clear whether the function $I_{[0,T]}$ has an integral form. Such an explicit expression of the rate function is important in view of applications.

1.3. *Results.* In the present paper, a complete proof of the upper large deviation bound (2) is given and an integral representation of the rate function $I_{[0,T]}$ is derived. The main arguments of the proof are now detailed.

For Markov processes considered in this paper, because of the partial homogeneity of the transitions, the local sample path large deviation principle is equivalent to the existence of a collection of convex nonnegative functions on $\mathbb{R}^N$,

$$L_\Lambda, \Lambda \subset \{1, \ldots, N\},$$

such that, for any $T > 0$ and for any linear function $\phi(s) = \phi(0) + vs$, $s \in [0, T]$, the sequence of scaled Markov processes $(Z_n(t), t \in [0, T])$ satisfies the inequalities

$$
\begin{aligned}
(6) \qquad w_{[0,T]}(\phi) &\stackrel{\text{def}}{=} \lim_{\delta \to 0} \lim_{\varepsilon \to 0} \liminf_{n \to \infty} \inf_{z : |z - \phi(0)| < \varepsilon} \frac{1}{n} \log \mathbb{P}_z(\|\phi - Z_n\|_\infty < \delta) \\
&\geq -T L_{\Lambda(\phi)}(v)
\end{aligned}
$$

and

$$
\begin{aligned}
(7) \qquad W_{[0,T]}(\phi) &\stackrel{\text{def}}{=} \lim_{\delta \to 0} \limsup_{n \to \infty} \sup_{z : |z - \phi(0)| < \delta} \frac{1}{n} \log \mathbb{P}_z(\|\phi - Z_n\|_\infty < \delta) \\
&\leq -T L_{\Lambda(\phi)}(v),
\end{aligned}
$$

where $\Lambda(\phi)$ is the set of all those $i \in \{1, \ldots, N\}$ for which $\phi_i(s) \neq 0$ for all $s \in (0, T)$. Inequalities (6) and (7) correspond to inequalities (3) and (4).

In order to get our main result, (7) will be replaced by a slightly different inequality, namely,

$$
\begin{aligned}
(8) \qquad \limsup_{n \to \infty} \sup_{z : |z - \phi(0)| < \delta} & \frac{1}{n} \log \mathbb{P}_z(|\phi(T) - Z_n(T)| < \delta \text{ and } \tau_{\Lambda,n} > T) \\
&\leq -T L_{\Lambda(\phi)}(v),
\end{aligned}
$$

where $\tau_{\Lambda,n}$ is the first time when the process $Z_n(t)$ hits the set $\bigcup_{i \in \Lambda}\{x : x_i = 0\}$.

When the local estimates (6) and (8) are verified and when the local rate function $L_{\{1,\ldots,N\}}$ is finite in a neighborhood of zero, it is proved that, for any absolutely continuous function $\phi : [0, T] \to \mathbb{R}_+^N$,

$$
(9) \qquad W_{[0,T]}(\phi) = w_{[0,T]}(\phi) = -I_{[0,T]}(\phi) = -\int_0^T L(\phi(t), \dot{\phi}(t)) \, dt.
$$

It is shown that this result implies the whole sample path large deviation principle when the general rough upper bound of [8] holds.

In the setting considered by Dupuis and Ellis [6], inequality (8) is proved in Section 3 by using the method of convergence parameters of corresponding local transform semigroups. This method was developed in [11] for partially homogeneous discrete-time Markov chains. In this way, our results complete the proof of the main result of [6].

The paper is organized as follows. Section 2 presents an overview of the main results. In Section 3, as an application, these results are used to establish the sample path large deviation principle for a general class of queueing networks. The proof of the local estimates (6) and (8) is given. Section 4 is devoted to the proof of (9) (this is the proof of Theorem 1 below). Using this relation and the general upper bound of [8], the whole sample path large deviation principle (Theorem 2) is derived in Section 5.

## 2. General results.

2.1. *Definitions and assumptions.* The following notation are used throughout this paper. For $x \in \mathbb{R}_+^N$, $\Lambda(x)$ is the set of those indices $i$ for which $x_i > 0$. For a subset $\Lambda \subset \{1, \ldots, N\}$:

(a) $x_\Lambda = (x_i, i \in \Lambda)$;
(b) $B_\Lambda$ is the set of all $x \in \mathbb{R}_+^N$ with $\Lambda(x) = \Lambda$.

It is assumed that the subsets $(\mathcal{E}_n, n \geq 1)$, the state spaces of scaled processes $(Z_n, n \geq 1)$, are dense in $\mathbb{R}_+^N$: for any $x \in \mathbb{R}_+^N$ there exists a sequence of points $x_n \in \mathcal{E}_n$ converging to $x$.

It is assumed that there is a collection of convex nonnegative functions $L_\Lambda$ on $\mathbb{R}^{\mathbb{N}}$ satisfying the following conditions.

$(A_1)$ For any $\Lambda \subset \{1, \ldots, N\}$ and $T > 0$, and for any function $\phi : [0, T] \to \mathbb{R}_+^N$ with a constant velocity $\dot\phi(t) = v$ and such that $\phi(t) \in B_\Lambda$ for $0 < t < T$, the following inequality holds:

$$\lim_{\delta \to 0} \lim_{\varepsilon \to 0} \liminf_{n \to \infty} \inf_{z \,:\, |z - \phi(0)| < \varepsilon} \frac{1}{n} \log \mathbb{P}_z(\|\phi - Z_n\|_\infty < \delta) \geq -T L_\Lambda(v).$$

$(A_2)$ The function $L_{\{1, \ldots, N\}}$ is finite in a neighborhood of zero.

$(A_3)$ For any $\Lambda \subset \{1, \ldots, N\}$ and $T > 0$, and for any $x, y \in B_\Lambda$, the following inequality holds:

$$\lim_{\delta \to 0} \limsup_{n \to \infty} \sup_{z \,:\, |z - x| < \delta} \frac{1}{n} \log \mathbb{P}_z(|Z_n(T) - y| < \delta \text{ and } \tau_{\Lambda, n} > T)$$
$$\leq -T L_\Lambda \left( \frac{y - x}{T} \right),$$

where $\tau_{\Lambda, n}$ is the hitting time of the set $\bigcup_{i \in \Lambda} \{x : x_i = 0\}$ by the process $(Z_n(t))$. In the next section we will see that these conditions are satisfied for a general class of queueing networks.

For a continuous piecewise linear function $\phi : [0, T] \to \mathbb{R}_+^N$, we define

$$J_{[0,T]}(\phi) = \int_0^T L(\phi(s), \dot\phi(s)) \, ds$$

with $L(x, v) = L_{\Lambda(x)}(v)$ for all $x \in \mathbb{R}_+^N$, $v \in \mathbb{R}^N$. The function $I_{[0,T]}$ is defined on $D([0, T], \mathbb{R}_+^N)$ by

$$I_{[0,T]}(\phi) = \lim_{\delta \to 0} \inf_\psi J_{[0,T]}(\psi),$$

the infimum being taken over all continuous piecewise linear functions $\psi$ on $[0, T]$ with $d_S(\phi, \psi) < \delta$, where $d_S(\cdot, \cdot)$ is the Skorohod metric on $D([0, T], \mathbb{R}_+^N)$. It is the lower semicontinuous regularization of the function $J_{[0,T]}$.

2.2. *The main theorems.* The central result of our paper is the following theorem.

THEOREM 1. *Under the assumptions* $(A_1)$, $(A_2)$ *and* $(A_3)$, *for any absolutely continuous function* $\phi : [0, T] \to \mathbb{R}_+^N$,

$$I_{[0,T]}(\phi) = \int_0^T L(\phi(t), \dot{\phi}(t)) \, dt = -W_{[0,T]}(\phi) = -w_{[0,T]}(\phi).$$

Recall that a mapping $\tilde{I}_{[0,T]} : D([0, T], \mathbb{R}_+^N) \to \mathbb{R}_+$ is a *good rate function*, if the following assertions hold:

(a) for any compact set $V \subset \mathbb{R}_+^N$ and for any $c > 0$ the set of all functions $\phi \in D([0, T], \mathbb{R}_+^N)$ with $\phi(0) \in V$ satisfying the inequality $\tilde{I}_{[0,T]}(\phi) \le c$ is compact in $D([0, T], \mathbb{R}_+^N)$;

(b) every function $\phi \in D([0, T], \mathbb{R}_+^N)$ with $\tilde{I}_{[0,T]}(\phi) < \infty$ is absolutely continuous.

A general upper large deviation bound with a good rate function was obtained by Dupuis, Ellis and Weiss [8]. The next theorem establishes that $I_{[0,T]}$ is a good rate function and that the sequence of Markov processes $(Z_n(t))$ satisfies the whole sample path large deviation principle when the general upper bound of [8] holds.

THEOREM 2. *Suppose that there is a good rate function* $\tilde{I}_{[0,T]}$ *satisfying the upper large deviation bound and let the hypotheses* $(A_1)$, $(A_2)$ *and* $(A_3)$ *be satisfied. Then, the sample path large deviation principle holds with the rate function* $I_{[0,T]}$ *and* $I_{[0,T]}$ *is also a good rate function.*

The main steps of our proofs are now briefly outlined. The proof of Theorem 1 begins by showing that for any absolutely continuous function $\phi : [0, T] \to \mathbb{R}_+^N$,

$$(10) \qquad I_{[0,T]}(\phi) \le \int_0^T L(\phi(t), \dot{\phi}(t)) \, dt.$$

To obtain this inequality, the classical approach consists in constructing for every $\varepsilon > 0$ a piecewise linear interpolation $\psi_\varepsilon$ of $\phi$ such that $\|\phi - \psi_\varepsilon\|_\infty < \varepsilon$ and

$$I_{[0,T]}(\psi_\varepsilon) = \int_0^T L(\psi_\varepsilon(t), \dot{\psi}_\varepsilon(t)) \, dt \le \int_0^T L(\phi(t), \dot{\phi}(t)) \, dt + \eta_\varepsilon,$$

where $\eta_\varepsilon \to 0$ as $\varepsilon \to 0$.

For Markov processes with a discontinuity in the transition mechanism along a hyperplane, such a construction was performed in Lemma 7.5.4

of [7] and in Lemma 4.9 of [3]. In some particular cases, when there is a nonnegative function $\ell$ on $\mathbb{R}^N$ such that $c_1 \ell(v) \le L(x, v) \le c_2 \ell(v)$ for all $x \in \mathbb{R}_+^N$ and for all $v \in \mathbb{R}^N$ with $v_i = 0$ for $i \notin \Lambda(x)$, this method can be extended to higher dimensions; see [2].

In our setting, such a construction does not seem possible: when $N \ge 3$ and when the trajectory $\{\phi(t), t \in [0, T]\}$ has a spiral form with an infinite number of linear segments on the boundary set $\bigcup_i \{x : x_i = 0\}$ converging to the center of the spiral $0 \in \mathbb{R}^N$, one can have $I_{[0,T]}(\psi) = +\infty$ for every piecewise linear interpolation $\psi$ of the function $\phi$.

Generally, a construction of the above piecewise linear interpolation is difficult and sometimes impossible in a neighborhood of some irregular points (in the above example, it is a center of the spiral). To avoid this difficulty, we slow down the velocity of piecewise linear interpolations in a neighborhood of irregular points. It is shown that, for any $\varepsilon > 0$, there is a piecewise linear interpolation $\psi_\varepsilon$ of $\phi$ and there is a strictly increasing continuous piecewise linear mapping $\theta_\varepsilon : [0, T] \to [0, \theta_\varepsilon(T)]$ with $\theta_\varepsilon(0) = 0$ and $\dot{\theta}_\varepsilon(t) \ge 1$ for almost all $t \in [0, T]$, such that

$$I_{[0,T]}(\psi_\varepsilon \circ \theta_\varepsilon^{-1}) \le \int_0^T L(\phi(t), \dot{\phi}(t)) \, dt + \varepsilon$$

and such that $\|\phi - \psi_\varepsilon\|_\infty$ and $\sup_{t \in [0,T]} |\theta_\varepsilon(t) - t|$ tend to 0 as $\varepsilon$ tends to 0. Since $\theta_\varepsilon(T) \ge T$, the resulting function $\psi_\varepsilon \circ \theta^{-1}$ is piecewise linear and continuous on $[0, T]$. The function $\phi$ being continuous, we obtain moreover that $\|\phi - \psi_\varepsilon \circ \theta_\varepsilon^{-1}\|_\infty$ converges to 0 as $\varepsilon$ tends to 0 and therefore, that (10) holds.

The next step is the proof of the inequality

$$W_{[0,T]}(\phi) \le -\int_0^T L(\phi(t), \dot{\phi}(t)) \, dt$$

for any absolutely continuous path $\phi$. To obtain this inequality, (8) is used.

The proof of the last inequality as well as the proof of the existence of $\psi_\varepsilon$ and $\theta_\varepsilon$ is performed by a careful induction with respect to $\Lambda \subset \{1, \ldots, N\}$ for $\phi = (\phi_1, \ldots, \phi_N) : [0, T] \to \mathbb{R}_+^N$ with $\phi_i(t) > 0$ for all $i \in \Lambda$ and for all $t \in [0, T]$.

Finally, with the lower large deviation bound of [6], we conclude that

$$-I_{[0,T]}(\phi) \le w_{[0,T]}(\phi) \le W_{[0,T]}(\phi) \le -\int_0^T L(\phi(t), \dot{\phi}(t)) \, dt \le -I_{[0,T]}(\phi).$$

This completes the proof of Theorem 1. Theorem 2 is proved classically with the results of Theorem 1.

## 3. Application: the large deviations of queueing networks.

In this section an application of our general results is presented to establish the sample path large deviation principle for Markov processes describing a general class of queueing networks.

For $x \in \mathbb{Z}_+^N$, we consider a continuous-time Markov process $(X(t, x))$ on $\mathbb{Z}_+^N$ generated by

$$\mathcal{L}f(y) = \sum_{y' \in \mathbb{Z}_+^N} q(y, y')(f(y') - f(y)), \qquad y \in \mathbb{Z}_+^N,$$

with $X(0, x) = x$. The transition intensities $q(y, y')$ of this process are assumed to satisfy the following conditions:

$(B_0)$ (*Finite range.*) There is $d > 0$ such that $q(y, y') = 0$ whenever $|y - y'| > d$.

$(B_1)$ (*Communication condition.*) There are $C > 0$ and $0 < \gamma < 1$ such that for any $y, y' \in \mathbb{Z}_+^N$, there exists a sequence $y_0 = y, y_1, \ldots, y_n = y' \in \mathbb{Z}_+^N$ with $n \leq C|y - y'|$ such that $q(y_{i-1}, y_i) \geq \gamma$ for all $i = 1, \ldots, n$.

$(B_2)$ (*Partial homogeneity.*) For every $\Lambda \subset \{1, \ldots, N\}$, there is a nonnegative measure $\mu_\Lambda$ on $\mathbb{Z}^N \setminus \{0\}$ such that

$$q(y, y') = \mu_\Lambda(y' - y)$$

for any $y \in \mathbb{Z}_+^N$ with $\Lambda(y) = \Lambda$, and for any $y' \in \mathbb{Z}_+^N$.

Recall that for $x \in \mathbb{R}_+^N$, $\Lambda(x)$ denotes the set of all those $i \in \{1, \ldots, N\}$ for which $x_i > 0$ and $B_\Lambda = \{x \in \mathbb{R}_+^N : \Lambda(x) = \Lambda\}$. For $\Lambda \subset \{1, \ldots, N\}$ and $x = (x_1, \ldots, x_N)$ we denote $x_\Lambda = (x_i; i \in \Lambda)$.

We prove that under the above assumptions, the sequence of scaled Markov processes

$$Z_n(t, z) = X(nt, nz)/n$$

satisfies the sample path large deviation principle in $D([0, T], \mathbb{R}_+^N)$ with a good rate function having an integral representation.

To prove the local large deviation estimates (6) and (8), we use the local Markov processes introduced in [6]. Roughly speaking, if the Markov process $(X(t))$ describes a queueing network with $N$ nodes, a local Markov process $(A_\Lambda(t), Y_\Lambda(t))$ on

$$\mathbb{Z}^\Lambda \times \mathbb{Z}_+^{\Lambda^c} = \{z \in \mathbb{Z}^N : z_i \geq 0 \text{ for all } i \in \Lambda^c\}$$

describes a modified queueing network with the same parameters as the original Markov process $(X(t))$, but without any boundary condition on the nodes $i \in \Lambda$: the queue lengths at the nodes $i \in \Lambda$ may be negative. Such a Markov process $(A_\Lambda(t), Y_\Lambda(t))$ is generated by

$$\mathcal{L}_\Lambda f(z) = \sum_{z' \in \mathbb{Z}^\Lambda \times \mathbb{Z}_+^{\Lambda^c}} q_\Lambda(z, z')(f(z') - f(z)),$$

where $q_\Lambda(z, z') = \mu_{\Lambda \cup \Lambda(z)}(z' - z)$. Throughout this section, we identify $(x_\Lambda, x_{\Lambda^c}) \in \mathbb{R}^\Lambda \times \mathbb{R}_+^{\Lambda^c}$ with $x = (x_1, \ldots, x_N)$.

The transition intensities $q_\Lambda(z, z')$ being invariant with respect to the translations on the first coordinate $z_\Lambda$, following the usual terminology, $(A_\Lambda(t), Y_\Lambda(t))$ is a *Markov-additive process* with additive part $A_\Lambda(t)$ on $\mathbb{Z}^\Lambda$ and with Markovian part $Y_\Lambda(t)$ on $\mathbb{Z}_+^{\Lambda^c}$. The Markovian part $Y_\Lambda(t)$ is a Markov chain on $\mathbb{Z}_+^{\Lambda^c}$. For $\Lambda = \{1, \ldots, N\}$, the Markovian part is empty and the local process $A_{\{1,\ldots,N\}}(t)$ is a random walk on $\mathbb{Z}^N$ with transition intensities $q_{\{1,\ldots,N\}}(z, z') = \mu_{\{1,\ldots,N\}}(z' - z)$.

The local estimates (6) and (8) are proved and the local rate function $L_\Lambda$ is expressed by using the method of convergence parameters of transform semigroups developed earlier in [11]. For a given $\alpha \in \mathbb{R}^\Lambda$, the *transform semigroup* $(\mathcal{P}_\Lambda^t(\alpha))$ of the Markov-additive process $(A_\Lambda(t), Y_\Lambda(t))$ is defined by

$$\mathcal{P}_\Lambda^t(\alpha)f(y) = \sum_{y' \in \mathbb{Z}^{\Lambda^c}} \mathcal{P}_\Lambda^t(\alpha; y, y')f(y') = \mathbb{E}_{(0,y)}(e^{\langle \alpha, A_\Lambda(t) \rangle} f(Y_\Lambda(t)))$$

for a nonnegative function $f : \mathbb{Z}_+^{\Lambda^c} \to \mathbb{R}$. $\mathbb{E}_{(0,y)}(\cdot)$ denotes here a conditional expectation given that $A_\Lambda(0) = 0$ and $Y_\Lambda(0) = y$. Under our assumptions, for all $\Lambda \subset \{1, \ldots, N\}$, $t > 0$ and $\alpha \in \mathbb{R}^\Lambda$, the quantities

$$\mathcal{P}_\Lambda^t(\alpha; y, y') = \mathbb{E}_{(0,y)}(e^{\langle \alpha, A_\Lambda(t) \rangle} \mathbb{1}_{\{Y_{\Lambda^c}(t) = y'\}}), \qquad y, y' \in \mathbb{Z}_+^{\Lambda^c},$$

are finite. Moreover, because of the communication condition $(B_1)$, the infinite matrices $\mathcal{P}_\Lambda^t(\alpha) = (\mathcal{P}_\Lambda^t(\alpha; y, y'); \, y, y' \in \mathbb{Z}_+^{\Lambda^c})$ are irreducible. Using the inequality

$$\mathcal{P}_\Lambda^{t+s+s'}(\alpha, y, y') \geq \mathcal{P}_\Lambda^s(\alpha, y, z)\mathcal{P}_\Lambda^t(\alpha, z, z')\mathcal{P}_\Lambda^{s'}(\alpha, z', y'),$$

this implies that the limit

$$\lambda_\Lambda(\alpha) = \limsup_{t \to \infty} \frac{1}{t} \log \mathcal{P}_\Lambda^t(\alpha; y, y')$$

does not depend on $y, y' \in \mathbb{Z}_+^{\Lambda^c}$. The quantity $\exp(-\lambda_\Lambda(\alpha))$ is called the *convergence parameter* of the semigroup $(\mathcal{P}_\Lambda^t(\alpha))$. For $\Lambda = \{1, \ldots, N\}$, clearly

$$\mathcal{P}_{\{1,\ldots,N\}}^t(\alpha) = \mathbb{E}_0(\exp\{\langle \alpha, A_{\{1,\ldots,N\}}(t) \rangle\})$$

$$= \exp\left( t \sum_{z \in \mathbb{Z} : z \neq 0} \mu_{\{1,\ldots,N\}}(z)(e^{\langle \alpha, z \rangle} - 1) \right)$$

and

$$\lambda_{\{1,\ldots,N\}}(\alpha) = \sum_{z \in \mathbb{Z} : z \neq 0} \mu_{\{1,\ldots,N\}}(z)(e^{\langle \alpha, z \rangle} - 1).$$

For $\Lambda \subset \{1, \ldots, N\}$, we define the function $L_\Lambda : \mathbb{R}^N \to \mathbb{R}$ by setting $L_\Lambda(v) = \lambda_\Lambda^*(v_\Lambda)$ where $\lambda_\Lambda^*$ is the convex conjugate of the function $\lambda_\Lambda$:

$$\lambda_\Lambda^*(v_\Lambda) = \sup_{\alpha \in \mathbb{R}^\Lambda} (\langle \alpha, v_\Lambda \rangle - \lambda_\Lambda(\alpha))$$

and we let $L(x, v) = L_\Lambda(v)$ if $\Lambda(x) = \Lambda$.

The main result of this section is the following theorem.

THEOREM 3. *Under the hypotheses* $(B_0)$, $(B_1)$ *and* $(B_2)$, *the sequence of scaled processes* $(Z_n(t))$ *satisfies the sample path large deviation principle with a good rate function*

$$I_{[0,T]}(\phi) = \begin{cases} \displaystyle\int_0^T L(\phi(t), \dot\phi(t))\,dt, & \text{if } \phi \text{ is absolutely continuous,} \\ +\infty, & \text{otherwise.} \end{cases}$$

The following lemmas prove the local estimates (6) and (8).

LEMMA 3.1. *For any* $\Lambda \subset \{1, \ldots, N\}$ *and* $T > 0$, *and for any linear path* $\phi(s) = \phi(0) + vt$ *with* $\phi(t) \in B_\Lambda$ *for* $s \in (0, T)$, (6) *holds.*

To prove this lemma it is sufficient to show that for any $\Lambda \subset \{1, \ldots, N\}$, $T > 0$ and $v \in \mathbb{R}^\Lambda$, the local Markov-additive process $(A_\Lambda(t), Y_\Lambda(t))$ satisfies the inequality

$$(11) \quad \lim_{\delta \to 0} \liminf_{n \to \infty} \frac{1}{n} \log \mathbb{P}_{(0,0)}\left( \sup_{t \in [0,nT]} |A_\Lambda(t) - tv| + |Y_\Lambda(t)| < \delta n \right) \geq -T\lambda_\Lambda^*(v)$$

(this is a consequence of Proposition 3.7 of [6]). In [11], this inequality was proved for discrete-time Markov-additive processes. For continuous-time Markov-additive processes the proof of (11) is quite similar. In the Appendix, we recall the main steps of this proof.

LEMMA 3.2. *For any* $\Lambda \subset \{1, \ldots, N\}$, $T > 0$ *and* $x, y \in B_\Lambda$, (8) *holds.*

PROOF. Remark that before the time $\tau_{\Lambda,n}$ when the process $Z_n(t) = X(nt)/n$ hits the set $\bigcup_{i \in \Lambda}\{x : x_i = 0\}$ for the first time, the transition intensities of the Markov process $X(t)$ are the same as those of the local Markov process $(A_\Lambda(t), Y_\Lambda(t))$. Hence, to prove (8) it is sufficient to show that the local Markov-additive process $(A_\Lambda(t), Y_\Lambda(t))$ satisfies the inequality

$$\begin{aligned}(12) \quad &\lim_{\delta \to 0} \limsup_{n \to \infty} \sup_{z \,:\, |z| < \delta n} \frac{1}{n} \log \mathbb{P}_z(|A_\Lambda(nT) - nTv| + |Y_\Lambda(nT)| < \delta n) \\ &\leq -T\lambda_\Lambda^*(v).\end{aligned}$$

For $\delta > 0$, $n \in \mathbb{N}$ and $v \in \mathbb{R}^\Lambda$, denote $E_{n\delta}(v) = \{|A_\Lambda(nT) - nTv| + |Y_\Lambda(nT)| < \delta n\}$. We will show that for any $\alpha \in \mathbb{R}^\Lambda$ such that $\lambda_\Lambda(\alpha) < +\infty$, and for any $\lambda > \lambda_\Lambda(\alpha)$,

$$(13) \quad \lim_{\delta \to 0} \limsup_{n \to \infty} \sup_{z \,:\, |z| < \delta n} \frac{1}{n} \log \mathbb{P}_z(E_{n\delta}(v)) \leq -T(\langle \alpha, v \rangle - \lambda)$$

from which (12) will follow.

Given $\alpha \in \mathbb{R}^\Lambda$ and $\lambda > \lambda_\Lambda(\alpha)$, we consider the function

$$f_\lambda(y) = \int_0^\infty \mathcal{P}_\Lambda^t(\alpha; y, 0) e^{-\lambda t}\, dt, \qquad y \in \mathbb{Z}_+^{\Lambda^c}.$$

According to the definition of $\lambda_\Lambda(\alpha)$, the above integral converges and for any $t > 0$,

$$(14) \qquad \mathcal{P}_\Lambda^t(\alpha) f_\lambda \leq e^{\lambda t} f_\lambda.$$

Furthermore, under the hypotheses $(B_0)$–$(B_2)$, there are $\mu$, $C_1$ and $C_2 > 0$ such that for any $y \in \mathbb{Z}_+^{\Lambda^c}$ there exists $n$ satisfying the inequalities $C_1|y| \leq n \leq C_2|y|$ and

$$\mathcal{P}^t(\alpha; y, 0) \geq \mathbb{P}_{(0,y)}(A_\Lambda(t) = 0, Y_\Lambda(t) = 0) \geq (\gamma t)^n e^{-\mu t}/n!$$

for all $t > 0$. This implies that

$$(15) \qquad f_\lambda(y) \geq \gamma^n/(\lambda + \mu)^{n+1}.$$

There exists moreover $m$ such that $C_1|y| \leq m \leq C_2|y|$ and

$$\mathcal{P}^t(\alpha; 0, y) \geq P_{(0,0)}(A_\Lambda(t) = 0, Y_\Lambda(t) = y) \geq (\gamma t)^m e^{-\mu t}/m!$$

for all $t > 0$. Hence, using (14) we obtain

$$(\gamma t)^m f_\lambda(y) e^{-\mu t}/m! \leq \mathcal{P}^t(\alpha; 0, y) f_\lambda(0) \leq e^{\lambda t} f_\lambda(0).$$

The last inequality with $t = m$ and (15) show that for any $\lambda > \lambda_\Lambda(\alpha)$ there is $c > 1$ such that

$$(16) \qquad c^{-|y|} \leq f_\lambda(y) \leq c^{|y|}$$

and hence, on the event $E_{n\delta}(v)$, the following inequality holds:

$$e^{\langle \alpha, A_\Lambda(nT)\rangle} f_\lambda(Y_\Lambda(nT)) \geq \exp(nT\langle \alpha, v\rangle - |\alpha|\delta n) c^{-\delta n}.$$

By Chebyshev's inequality, this implies that

$$\mathbb{P}_z(E_{n\delta}(v)) \leq c^{\delta n} \exp(|\alpha|\delta n - nT\langle \alpha, v\rangle) \mathbb{E}_z(e^{\langle \alpha, A_\Lambda(nT)\rangle} f_\lambda(Y_\Lambda(nT))).$$

Moreover, using (14) it follows that

$$\mathbb{E}_z(e^{\langle \alpha, A_\Lambda(nT)\rangle} f_\lambda(Y_\Lambda(nT))) = e^{\langle \alpha, z_\Lambda\rangle} \mathcal{P}^{nT}(\alpha) f_\lambda(z_{\Lambda^c}) \leq e^{\langle \alpha, z_\Lambda\rangle} e^{\lambda nT} f_\lambda(z_{\Lambda^c})$$

and consequently, using again (16) we obtain

$$\sup_{z\,:\,|z| < \delta n} \frac{1}{n} \log \mathbb{P}_z(E_{n\delta}(v)) \leq 2\delta \log c + 2|\alpha|\delta - T\langle \alpha, v\rangle + \lambda T.$$

Letting $n \to \infty$ and $\delta \to 0$ in the last equality, (13) follows. Moreover, letting $\lambda \to \lambda_\Lambda(\alpha)$ in (13), it follows that for all $\alpha \in \mathrm{dom}(\lambda_\Lambda) = \{\alpha : \lambda_\Lambda(\alpha) < +\infty\}$,

$$\lim_{\delta \to 0} \limsup_{n \to \infty} \sup_{z \,:\, |z| < \delta n} \frac{1}{n} \log \mathbb{P}_z(E_{n\delta}(v)) \leq -T(\langle \alpha, v \rangle - \lambda_\Lambda(\alpha))$$

and hence,

$$\lim_{\delta \to 0} \limsup_{n \to \infty} \sup_{z \,:\, |z| < \delta n} \frac{1}{n} \log \mathbb{P}(E_{n\delta}(v)) \leq -T \sup_{\alpha \in \mathrm{dom}(\lambda_\Lambda)} (\langle \alpha, v \rangle - \lambda_\Lambda(\alpha)).$$

The last inequality proves (12) because

$$\lambda_\Lambda^*(v) = \sup_{\alpha \in \mathbb{R}^\Lambda} (\langle \alpha, v \rangle - \lambda_\Lambda(\alpha)) = \sup_{\alpha \in \mathrm{dom}(\lambda_\Lambda)} (\langle \alpha, v \rangle - \lambda_\Lambda(\alpha))$$

(see [12], Corollary 12.2.2 of Theorem 12.2).   □

PROOF OF THEOREM 3.   We are ready now to prove Theorem 3. For this, it is sufficient to show that the hypotheses of Theorem 2 are satisfied.

Conditions $(A_1)$ and $(A_3)$ are satisfied because of Lemmas 3.1 and 3.2. Moreover, under the hypotheses $(B_0)$ and $(B_1)$, the convex conjugate of the function

$$\lambda_{\{1,\dots,N\}}(\alpha) = \sum_{z \in \mathbb{Z} \,:\, z \neq 0} \mu_{\{1,\dots,N\}}(z)(e^{\langle \alpha, z \rangle} - 1)$$

is finite in a neighborhood of zero and consequently, the condition $(A_2)$ is also satisfied. Finally, under our hypotheses, the general upper large deviation bound of [8] holds and hence, Theorem 2 can be applied.   □

## 4. Proof of Theorem 1.

Let $D([a, b], \mathbb{R}_+^N)$ be the set of all functions $\phi : [a, b] \to \mathbb{R}_+^N$ which are right-continuous and have left limits. It is convenient to introduce the functions $I_{[a,b]}(\cdot)$, $w_{[a,b]}(\cdot)$ and $W_{[a,b]}(\cdot)$ on $D([a, b], \mathbb{R}_+^N)$ for every interval $[a, b] \subset \mathbb{R}_+$. For $\phi \in D([a, b], \mathbb{R}_+^N)$, the expressions $w_{[a,b]}(\phi)$ and $W_{[a,b]}(\phi)$ are generalized as follows:

$$w_{[a,b]}(\phi) \stackrel{\mathrm{def}}{=} \lim_{\delta \to 0} \lim_{\varepsilon \to 0} \liminf_{n \to \infty} \inf_{z \,:\, |z - \phi(t)| < \varepsilon} \frac{1}{n} \log \mathbb{P}_{a,z}\left( \sup_{s \in [a,b]} |\phi(s) - Z_n(s)| < \delta \right)$$

and

$$W_{[a,b]}(\phi) \stackrel{\mathrm{def}}{=} \lim_{\delta \to 0} \limsup_{n \to \infty} \sup_{z \,:\, |z - \phi(t)| < \delta} \frac{1}{n} \log \mathbb{P}_{a,z}\left( \sup_{s \in [a,b]} |\phi(s) - Z_n(s)| < \delta \right),$$

where $\mathbb{P}_{a,z}$ is a conditional probability given that $Z_n(a) = z \in \mathcal{E}_n$.

Recall that a continuous function $\phi\colon [a,b] \to \mathbb{R}_+^N$ is called piecewise linear if there is $n \geq 1$ and there are $a = t_0 \leq t_1 \leq \cdots \leq t_n = b$ such that for all $t \in [t_{i-1}, t_i]$, $i = 1, \ldots, n$,

$$\phi(t) = \phi(t_{i-1}) + (t - t_{i-1}) \frac{\phi(t_i) - \phi(t_{i-1})}{t_i - t_{i-1}}.$$

For a continuous piecewise linear function $\psi\colon [a,b] \to \mathbb{R}_+^N$ we let

$$J_{[a,b]}(\psi) = \int_a^b L(\psi(t), \dot{\psi}(t)) \, dt.$$

The function $I_{[a,b]}$ is defined by

$$I_{[a,b]}(\phi) = \lim_{\delta \to 0} \inf_{\psi \,:\, d_S(\phi,\psi) < \delta} J_{[a,b]}(\psi),$$

the infimum being taken over all continuous piecewise linear $\psi\colon [a,b] \to \mathbb{R}_+^N$ with $d_S(\phi, \psi) < \delta$ where $d_S(\cdot, \cdot)$ is the Skorohod metric on $D([a,b], \mathbb{R}_+^N)$.

We begin the proof of Theorem 1 with the following proposition.

PROPOSITION 4.1.  *Under the hypotheses* $(A_1)$–$(A_3)$, *for any absolutely continuous function* $\phi\colon [a,b] \to \mathbb{R}_+^N$,

$$(17) \qquad I_{[a,b]}(\phi) \leq \int_a^b L(\phi(t), \dot{\phi}(t)) \, dt.$$

Recall that a piecewise linear function $\psi$ is called a piecewise linear interpolation of the function $\phi \in D([a,b], \mathbb{R}_+^N)$ if there is $n \geq 1$ and there are $a = t_0 \leq t_1 \leq \cdots \leq t_n = b$ such that for all $t \in [t_{i-1}, t_i]$, $i = 1, \ldots, n$,

$$\psi(t) = \phi(t_{i-1}) + (t - t_{i-1}) \frac{\phi(t_i) - \phi(t_{i-1})}{t_i - t_{i-1}}.$$

To obtain (17), we show that for any $\varepsilon > 0$ there is a piecewise linear interpolation $\psi_\varepsilon$ of $\phi$ and there is a strictly increasing continuous piecewise linear function $\theta_\varepsilon\colon [a,b] \to \mathbb{R}_+$ with $\theta_\varepsilon(a) = a$ and $\theta_\varepsilon(b) \geq b$, such that

$$I_{[a,\theta_\varepsilon(b)]}(\psi_\varepsilon \circ \theta_\varepsilon^{-1}) \leq \int_a^b L(\phi(t), \dot{\phi}(t)) \, dt + \varepsilon,$$

and such that $\sup_{t \in [a,b]} |t - \theta_\varepsilon(t)| \to 0$ and $\|\phi - \psi_\varepsilon\| \to 0$ when $\varepsilon \to 0$. Then, the function $\phi$ being continuous,

$$
\begin{aligned}
(18) \quad & \|\psi_\varepsilon \circ \theta_\varepsilon^{-1} - \phi\|_\infty \\
&= \sup_{t \in [a,b]} \|\psi_\varepsilon \circ \theta_\varepsilon^{-1}(t) - \phi(t)\| \\
&= \sup_{t \in [a,\theta_\varepsilon^{-1}(b)]} \|\psi_\varepsilon(t) - \phi \circ \theta_\varepsilon(t)\| \\
&\leq \|\psi_\varepsilon - \phi\|_\infty + \sup_{t \in [a,\theta_\varepsilon^{-1}(b)]} \|\phi(t) - \phi \circ \theta_\varepsilon(t)\| \to 0 \qquad \text{as } \varepsilon \to 0
\end{aligned}
$$

and hence, (17) will follow.

For our purpose, it is convenient to introduce a new function $G_{[a,b]}$ by letting

$$(19) \qquad G_{[a,b]}(\phi) = \liminf_{\delta \to 0} \inf_{\psi, \theta} I_{[a,\theta(b)]}(\psi \circ \theta^{-1}),$$

where the infimum is taken over all piecewise linear interpolations $\psi$ of $\phi$ such that $\|\phi - \psi\|_\infty < \delta$ and over all continuous piecewise linear functions $\theta \colon [a, b] \to \mathbb{R}$ such that $\theta(a) = a$, $\sup_{t \in [a,b]} |t - \theta(t)| < \delta$ and $\dot{\theta}(t) \geq 1$ for almost all $t \in [a, b]$. To prove Proposition 4.1 we will use the following properties of the function $G_{[a,b]}$.

LEMMA 4.1.   *For any continuous function* $\phi \colon [a, b] \to \mathbb{R}_+^N$, *and for any* $c \in [a, b]$,

$$(20) \qquad I_{[a,b]}(\phi) \leq G_{[a,b]}(\phi) \leq G_{[a,c]}(\phi) + G_{[c,b]}(\phi).$$

PROOF.   The first inequality of (20) follows from (18). The second inequality holds because

$$G_{[a,c]}(\phi) + G_{[c,b]}(\phi) = \liminf_{\delta \to 0} \inf_{\psi, \theta} I_{[a,\theta(b)]}(\psi \circ \theta^{-1}),$$

where the infimum is taken over all $\psi$ and $\theta$ satisfying the same condition as in (19) but with $\psi(c) = \phi(c)$.   □

To prove the next property of the function $G_{[a,b]}$ we need the following lemma.

LEMMA 4.2.   *For any* $\Lambda \subset \{1, \ldots, N\}$ *and for any* $v \in \mathbb{R}^N$ *with* $v_{\Lambda^c} = 0$,

$$(21) \qquad L_\Lambda(v) \leq L_{\{1,\ldots,N\}}(v).$$

PROOF.   Let $\Lambda \subset \{1, \ldots, N\}$ and let $v \in \mathbb{R}^N$ be such that $v_{\Lambda^c} = 0$. Consider $x \in B_\Lambda$ and $T > 0$ such that $\phi(t) = x + vt \in B_\Lambda$ for all $t \in [0, T]$. Then because of assumptions $(A_1)$ and $(A_3)$, the following relations hold:

$$w_{[0,T]}(\phi) = W_{[0,T]}(\phi) = -T L_\Lambda(v).$$

Similarly for $\phi_n(t) = \phi(t) + z/n$ with $z = (1, \ldots, 1)$,

$$w_{[0,T]}(\phi_n) = W_{[0,T]}(\phi_n) = -T L_{\{1,\ldots,N\}}(v).$$

The mapping $\phi \to W_{[0,T]}(\phi)$ being upper semicontinuous, this proves (21).   □

LEMMA 4.3.   *For any continuous function* $\phi \colon [a, b] \to \mathbb{R}_+^N$,

$$(22) \qquad G_{[a,b]}(\phi) \leq \lim_{\varepsilon \to 0^+} G_{[a+\varepsilon, b-\varepsilon]}(\phi).$$

PROOF.  By definition, for any continuous function $\phi = (\phi_1, \dots, \phi_N) : [a, b] \to \mathbb{R}_+^N$,

$$G_{[a,b]}(\phi) = \lim_{\delta \to 0^+} \inf_{\{t_i\}, \{\theta_i\}} \sum_{i=j}^n \theta_j (t_j - t_{j-1}) L_{\Lambda_j} \left( \frac{\phi(t_j) - \phi(t_{j-1})}{\theta_j (t_j - t_{j-1})} \right),$$

where for every $j = 1, \dots, n$, $\Lambda_j$ is the set of all those $i \in \{1, \dots, N\}$ for which

$$\phi_i(t_{j-1}) + (t - t_{j-1}) \frac{\phi_i(t_j) - \phi_i(t_{j-1})}{t_j - t_{j-1}} > 0 \qquad \text{for } t_{j-1} < t < t_j,$$

and the infimum is taken over all partitions $a = t_0 < t_1 < \cdots < t_n = b$ with $\max_i(t_i - t_{i-1}) < \delta$ and over all real numbers $\theta_i \geq 1$, $i = 1, \dots, n$, such that

$$\sum_{i=1}^n \theta_i(t_i - t_{i-1}) \leq b - a + \theta.$$

Letting $t_1 - t_0 = \varepsilon$ and $t_n - t_{n-1} = \varepsilon'$ and using Lemma 4.2, it follows therefore that

$$\begin{aligned}
G_{[a,b]}(\phi) \leq \lim_{\delta \to 0} \inf_{\varepsilon, \varepsilon', \theta, \theta'} &\theta \varepsilon L_{\{1, \dots, N\}} \left( \frac{\phi(a + \varepsilon) - \phi(a)}{\theta \varepsilon} \right) + G_{[a+\varepsilon, b-\varepsilon']}(\phi) \\
&+ \theta' \varepsilon' L_{\{1, \dots, N\}} \left( \frac{\phi(b) - \phi(b - \varepsilon')}{\theta' \varepsilon'} \right),
\end{aligned} \tag{23}$$

where the infimum is taken over all $\varepsilon, \varepsilon' > 0$, $\theta \geq 1$ and $\theta' \geq 1$ with $\varepsilon \theta + \varepsilon' \theta' < \delta$. Recall that by assumption $(A_2)$, the function $L_{\{1, \dots, N\}}$ is finite in a neighborhood of $0 \in \mathbb{R}^N$. Being convex it is therefore bounded in a neighborhood of $0 \in \mathbb{R}^N$ and hence, there are two real numbers $r > 0$ and $c > 0$ such that $L_{\{1, \dots, N\}}(v) \leq c$ for all $v \in \mathbb{R}^N$ with $|v| \leq r$. Without any restriction of generality we suppose that $r < 1$ and $c > 1$. For given $\delta > 0$, let us choose $0 < \varepsilon_\delta < \delta/(2c)$ such that for $0 < \varepsilon \leq \varepsilon_\delta$ $|\phi(a + \varepsilon) - \phi(a)| < r\delta/(2c)$ and let $\theta = \max\{1, |\phi(a + \varepsilon) - \phi(a)|/(r\varepsilon)\}$. Then $|\phi(a + \varepsilon) - \phi(a)|/(\varepsilon \theta) \leq r$ and hence,

$$\theta \varepsilon L_{\{1, \dots, N\}} \left( \frac{\phi(a + \varepsilon) - \phi(a)}{\theta \varepsilon} \right) \leq \theta \varepsilon c = c \max\{\varepsilon, |\phi(a + \varepsilon) - \phi(a)|/r\} \leq \delta/2.$$

The same arguments show that there are $\varepsilon_\delta' > 0$ and $\theta'(\varepsilon) \geq 1$ such that for $0 < \varepsilon' < \varepsilon_\delta'$

$$\begin{aligned}
\theta' \varepsilon' L_{\{1, \dots, N\}} &\left( \frac{\phi(b) - \phi(b - \varepsilon')}{\theta' \varepsilon'} \right) \\
&\leq \theta' \varepsilon' c = c \max\{\varepsilon', |\phi(b) - \phi(b - \varepsilon')|/r\} \leq \delta/2.
\end{aligned}$$

For such $\varepsilon, \varepsilon', \theta$ and $\theta'$, we have $\theta \varepsilon + \theta' \varepsilon' < \theta$ and hence, using (23) we obtain (22).  $\square$

LEMMA 4.4. *For any $\Lambda \subset \{1, \ldots, N\}$ and for any $\phi = (\phi_1, \ldots, \phi_N) \in D([a,b], \mathbb{R}_+^N)$ such that $\phi(a), \phi(b) \in B_\Lambda$ and $\phi_i(t) > 0$ for all $i \in \Lambda$ and for all $t \in [a,b]$, the following inequality holds:*

$$(24) \qquad (b-a)L_\Lambda\left(\frac{\phi(b) - \phi(a)}{b-a}\right) \leq I_{[a,b]}(\phi).$$

PROOF. Indeed, let $x, y \in B_\Lambda$ and let $\mathcal{O}_\delta$ be the set of all $\phi \in D([a,b], \mathbb{R}_+^N)$ with $|\phi(b) - y| < \delta$ and such that $\phi_i(t) > 0$ for all $i \in \Lambda$ and for all $t \in [a,b]$. Then because of assumption $(A_3)$,

$$(25) \quad \lim_{\delta \to 0} \limsup_{n \to \infty} \sup_{y\,:\,|y-x|<\delta} \frac{1}{n} \log \mathbb{P}_{a,y}(Z_n(\cdot) \in \mathcal{O}_\delta) \leq -(b-a)L_\Lambda\left(\frac{y-x}{b-a}\right).$$

Moreover, recall that under the hypotheses $(A_1)$, the rate function $I_{[0,b-a]}$ satisfies the lower large deviation bound (1) with $T = b - a$. The set $\mathcal{O}_\delta$ being open, using the Markov property we obtain

$$-I_{[a,b]}(\phi) \leq \lim_{\varepsilon \to 0} \liminf_{n \to \infty} \inf_{y\,:\,|y-x|<\varepsilon} \frac{1}{n} \log \mathbb{P}_{a,y}(Z_n(\cdot) \in \mathcal{O}_\delta)$$

for any $\phi \in \mathcal{O}_\delta$ with $\phi(a) = x$ and $\phi(b) = y$. Letting in the last inequality $\delta \to 0$ and using (25), (24) follows. $\square$

The last lemma combined with Lemma 4.1 implies the following property of the function $G_{[a,b]}$.

LEMMA 4.5. *For any $\Lambda \subset \{1, \ldots, N\}$ and for any $\phi = (\phi_1, \ldots, \phi_N) \in D([a,b], \mathbb{R}_+^N)$ such that $\phi(a), \phi(b) \in B_\Lambda$ and $\phi_i(t) > 0$ for all $i \in \Lambda$ and for all $t \in [a,b]$, the following inequality holds:*

$$(26) \qquad (b-a)L_\Lambda\left(\frac{\phi(b) - \phi(a)}{b-a}\right) \leq G_{[a,b]}(\phi).$$

PROOF OF PROPOSITION 4.1. We are ready now to prove Proposition 4.1. Because of Lemma 4.1, it is sufficient to show that for any absolutely continuous function $\phi = (\phi_1, \ldots, \phi_N) : [a,b] \to \mathbb{R}_+^N$, the following inequality holds:

$$(27) \qquad G_{[a,b]}(\phi) \leq \int_a^b L(\phi(t), \dot{\phi}(t))\, dt.$$

Suppose first that $\phi_i(t) > 0$ for all $i = 1, \ldots, N$ and for all $t \in [a,b]$; then

$$\int_a^b L(\phi(t), \dot{\phi}(t))\, dt = \int_a^b L_{\{1, \ldots, N\}}(\dot{\phi}(t))\, dt.$$

Moreover, for any piecewise linear interpolation $\psi = (\psi_1, \ldots, \psi_N)$ of the function $\phi$, we have also $\psi_i(t) > 0$ for all $i = 1, \ldots, N$ and for all $t \in [a, b]$ which implies that

$$I_{[a,b]}(\psi) = \int_a^b L(\psi(t), \dot{\psi}(t)) \, dt = \int_a^b L_{\{1,\ldots,N\}}(\dot{\psi}(t)) \, dt.$$

The function $L_{\{1,\ldots,N\}}(\cdot)$ being convex, this implies that

$$I_{[a,b]}(\psi) = \int_a^b L_{\{1,\ldots,N\}}(\dot{\psi}(t)) \, dt \leq \int_a^b L_{\{1,\ldots,N\}}(\dot{\phi}(t)) \, dt = \int_a^b L(\phi(t), \dot{\phi}(t)) \, dt$$

and hence, (27) holds.

To prove (27) in the general case, let us consider for every $\Lambda \subset \{1, \ldots, N\}$, the set $\Phi_\Lambda$ of all absolutely continuous functions $\phi = (\phi_1, \ldots, \phi_N) : [a, b] \to \mathbb{R}_+^N$ with arbitrary $a < b$ such that $\phi_i(t) > 0$ for all $t \in [a, b]$ and for all $i \in \Lambda$. We prove (27) by induction with respect to $\Lambda$ for $\phi \in \Phi_\Lambda$. Remark that for all $\phi \in \Phi_\Lambda$ with $\Lambda = \{1, \ldots, N\}$, this inequality is already verified.

Suppose that (27) is already verified for all $\phi \in \Phi_{\Lambda'}$ with $\Lambda' \subset \{1, \ldots, N\}$ such that $\Lambda \subset \Lambda' \neq \Lambda$ and let us verify this inequality for $\phi \in \Phi_\Lambda$, $\phi = (\phi_1, \ldots, \phi_N) : [a, b] \to \mathbb{R}_+^N$. Because of our assumption, for $a \leq t \leq t' \leq b$, the inequality

$$(28) \qquad\qquad G_{[t,t']}(\phi) \leq \int_t^{t'} L(\phi(s), \dot{\phi}(s)) \, ds$$

is already verified if there exists $i \in \{1, \ldots, N\} \setminus \Lambda$ such that $\phi_i(s) > 0$ for all $s \in [t, t']$.

Consider first the case when $\phi(t) = (\phi_1(t), \ldots, \phi_N(t)) \notin B_\Lambda$ for all $t \in [a, b]$. Then there is $\varepsilon > 0$ such that

$$(29) \qquad\qquad \sum_{i \notin \Lambda} \phi_i(t) > \varepsilon \qquad \text{for all } t \in [a, b]$$

and there is $\sigma > 0$ such that for all $t, s \in [a, b]$ satisfying the inequality $|t - s| < \sigma$, the inequality $\sum_i |\phi_i(t) - \phi_i(s)| < \varepsilon / N$ holds. Consider an increasing sequence $a = t_0 < t_1 < \cdots < t_n = b$ with $\sup_l |t_{l+1} - t_l| < \sigma$. If for $t \in [t_{l-1}, t_l]$, $\phi_i(t) = 0$ for some $i \in \{1, \ldots, N\} \setminus \Lambda$, then because of (29) there is $j \in \{1, \ldots, N\} \setminus \Lambda$ such that $\phi_j(t) > \varepsilon / (N - 1)$ and consequently, for any $s \in [t_{l-1}, t_l]$, the following inequality holds:

$$\phi_j(s) > \frac{\varepsilon}{N - 1} - |\phi_j(t) - \phi_j(s)| > \frac{\varepsilon}{N - 1} - \frac{\varepsilon}{N} > 0.$$

This proves that for any $l = 1, \ldots, n$, there is $j_l \in \{1, \ldots, N\} \setminus \Lambda$ such that $\phi_{j_l}(s) > 0$ for all $s \in [t_{l-1}, t_l]$ and hence, using (28) with $[t, t'] = [t_{l-1}, t_l]$ for each $l = 1, \ldots, n$, we obtain

$$G_{[a,t_1]}(\phi) + G_{[t_1,t_2]}(\phi) + \cdots + G_{[t_{n-1},b]}(\phi) \leq \int_a^b L(\phi(t), \dot{\phi}(t)) \, dt.$$

The last inequality and Lemma 4.1 imply (27).

Consider now an arbitrary function $\phi \in \Phi_\Lambda$, $\phi = (\phi_1, \ldots, \phi_N) : [a, b] \to \mathbb{R}_+^N$. Remark that for such a function $\phi$, (28) is already verified if $\phi(s) \notin B_\Lambda$ for all $s \in [t, t']$. The function $\phi$ being continuous, the set

$$\Delta = \{t \in (a, b) : \phi(t) \notin B_\Lambda\} = \bigcup_{i \in \Lambda} \{t \in (a, b) : \phi_i(t) > 0\}$$

is open and hence, it is a union of a countable family of open disjoint intervals $(t_k, t_k')$, $k \in \mathbb{N}$. For any $k \in \mathbb{N}$, and for any $0 < \sigma < (t_k' - t_k)/2$, the inequality

$$G_{[t_k+\sigma, t_k'-\sigma]}(\phi) \geq \int_{t_k+\sigma}^{t_k'-\sigma} L(\phi(s), \dot\phi(s)) \, ds$$

is already verified and hence, using Lemma 4.3 it follows that

$$G_{[t_k, t_k']}(\phi) \leq \lim_{\sigma \to 0} G_{[t_k+\sigma, t_k'-\sigma]}(\phi) \leq \int_{t_k}^{t_k'} L(\phi(s), \dot\phi(s)) \, ds.$$

According to the definition of the function $G_{[t_k, t_k']}(\phi)$ this implies that for given $\epsilon > 0$, there is a piecewise linear interpolation $\psi_k$ of the function $\phi : [t_k, t_k'] \to \mathbb{R}_+^N$ and a continuous piecewise linear function $\theta_k : [t_k, t_k'] \to \mathbb{R}$ with $\theta_k(t_k) = t_k$ and $\dot\theta_k(t) \geq 1$ for almost all $t \in [t_k, t_k']$, such that

$$\sup_{t \in [t_k, t_k']} |\phi(t) - \psi_k(t)| < \epsilon, \qquad \sup_{t \in [t_k, t_k']} |\theta_k(t) - t| < \epsilon/2^k$$

and

$$(30) \qquad I_{[\theta(t_k), \theta(t_k')]}(\psi_k \circ \theta_k^{-1}) \leq \int_{t_k}^{t_k'} L(\phi(s), \dot\phi(s)) \, ds + \varepsilon/2^k.$$

Moreover, Lemma 4.5 shows that for any $k \in \mathbb{N}$, for which $\phi(t_k), \phi(t_k') \in B_\Lambda$, the following inequality holds:

$$(31) \qquad (t_k' - t_k) L_\Lambda \left( \frac{\phi(t_k') - \phi(t_k)}{t_k' - t_k} \right) \leq G_{[t_k, t_k']}(\phi).$$

Given $\epsilon > 0$ let us choose $n_\epsilon$ such that

$$(32) \qquad \sum_{k \geq n_\epsilon} (t_k' - t_k) < \epsilon/2$$

and such that $n_\epsilon \geq k$ if $t_k = a$ or $t_k' = b$. Then for all $k > n_\varepsilon$, (31) holds and consequently,

$$(t_k' - t_k) L_\Lambda \left( \frac{\phi(t_k') - \phi(t_k)}{t_k' - t_k} \right) \leq \int_{t_k}^{t_k'} L(\phi(s), \dot\phi(s)) \, ds.$$

For the function

$$\phi_\varepsilon(t) = \begin{cases} \phi(t_k) + (t - t_k)(\phi(t'_k) - \phi(t_k))/(t'_k - t_k), \\ \qquad\qquad \text{for } t \in (t_k, t'_k), \ k > n_\varepsilon, \\ \phi(t), \qquad \text{for } t \in [a, b] \setminus \bigcup_{k > n_\varepsilon} (t_k, t'_k), \end{cases}$$

the above inequality implies that

$$(33) \quad \int_{[a,b] \setminus \bigcup_{k=1}^{n_\varepsilon}(t_k, t'_k)} L(\phi_\varepsilon(s), \dot\phi_\varepsilon(s)) \, ds \leq \int_{[a,b] \setminus \bigcup_{k=1}^{n_\varepsilon}(t_k, t'_k)} L(\phi(s), \dot\phi(s)) \, ds.$$

The set $(a, b) \setminus \bigcup_{k=1}^{n_\varepsilon} [t_k, t'_k]$ is a union of a finite number of disjoint open intervals $(s_i, s'_i)$, $i = 1, \ldots, m$. By construction, $\phi(s_i), \phi(s'_i) \in B_\Lambda$ for any $i = 1, \ldots, m$. For every $i = 1, \ldots, m$, we define a partition $s_i = s_{i0} < s_{i1} < \cdots < s_{ik_i} = s'_i$ by induction: if $s_{ij}$ is already defined:

(a) we let $s_{ij+1} = s'_i$ and $k_i = j + 1$, if $s'_i < s_{ij} + \epsilon$;

(b) otherwise, (32) shows that there is $s_{ij} + \epsilon/2 < s < s_{ij} + \epsilon$ such that $\phi(s) \in B_\Lambda$ and we let $s_{ij+1} = s$.

Then the piecewise linear function

$$\xi(t) = \phi(s_{ij-1}) + (t - s_{ij-1}) \frac{\phi(s_{ij}) - \phi(s_{ij-1})}{s_{ij} - s_{ij-1}},$$

$$t \in [s_{ij-1}, s_{ij}], \ j = 1, \ldots, k_i,$$

satisfies the following relations:

$$\int_{s_i}^{s'_i} L(\xi(t), \dot\xi(t)) \, dt = \int_{s_i}^{s'_i} L_\Lambda(\dot\xi(t)) \, dt$$

$$\leq \int_{s_i}^{s'_i} L_\Lambda(\dot\phi_\varepsilon(t)) \, dt$$

$$= \int_{s_i}^{s'_i} L(\phi_\varepsilon(t), \dot\phi_\varepsilon(t)) \, dt.$$

The first relation holds here because $\xi(t) \in B_\Lambda$ for all $t \in [s_{ij-1}, s_{ij}]$; the second relation is verified because by construction, $\xi$ is a piecewise linear interpolation of $\phi_\varepsilon$ and because the function $L_\Lambda(\cdot)$ is convex. Finally, the last identity is verified because $\phi_\varepsilon(t) \in B_\Lambda$ for all $t \in [s_{ij-1}, s_{ij}]$. Using (33) we conclude that

$$(34) \quad \int_{[a,b] \setminus \bigcup_{k=1}^{n_\epsilon}(t_k, t'_k)} L(\xi(s), \dot\xi(s)) \, ds \leq \int_{[a,b] \setminus \bigcup_{k=1}^{n_\epsilon}(t_k, t'_k)} L(\phi(s), \dot\phi(s)) \, ds.$$

Define now a piecewise linear interpolation $\psi_\varepsilon$ of the function $\phi$ on the whole interval $[a, b]$ by

$$\psi_\varepsilon(t) = \begin{cases} \psi_k(t), & \text{for } t \in (t_k, t_k'), \; k \le n_\varepsilon, \\ \xi(t), & \text{for } t \in [a, b] \setminus \bigcup_{k \le n_\varepsilon} (t_k, t_k'), \end{cases}$$

and let $\theta_\varepsilon \colon [a, b] \to \mathbb{R}$ be a continuous piecewise linear function with $\theta_\varepsilon(a) = a$ and

$$\dot{\theta}_\varepsilon(t) = \begin{cases} \dot{\theta}_k(t), & \text{for } t \in (t_k, t_k'), \; k \le n_\varepsilon, \\ 1, & \text{for } t \in (a, b) \setminus \bigcup_{k \le n_\varepsilon} [t_k, t_k']. \end{cases}$$

Then clearly, $\sup_{t \in [a, b]} |\theta_\varepsilon(t) - t| < \varepsilon$ and $\|\phi - \psi_\varepsilon\|_\infty \to 0$ when $\varepsilon \to 0$. Moreover, (30) and (34) imply that

$$I_{[\theta_\varepsilon(a), \theta_\varepsilon(b)]}(\psi_\varepsilon \circ \theta_\varepsilon^{-1}) = \int_{[a,b] \setminus \bigcup_{k=1}^{n_\varepsilon} (t_k, t_k')} L(\xi(s), \dot{\xi}(s)) \, ds + \sum_{i=1}^{n_\varepsilon} I_{[\theta_k(t_k), \theta_k(t_k')]}(\psi_k)$$

$$\le \int_a^b L(\phi(t), \dot{\phi}(t)) \, dt + \varepsilon$$

and hence, letting $\varepsilon \to 0$, (27) follows. Proposition 4.1 is therefore proved. $\square$

PROPOSITION 4.2. *Under the hypotheses* $(A_1)$ *and* $(A_3)$, *for any absolutely continuous function* $\phi \colon [a, b] \to \mathbb{R}_+^N$,

$$(35) \qquad W_{[a,b]}(\phi) \le - \int_a^b L(\phi(t), \dot{\phi}(t)) \, dt.$$

PROOF. Let $t \to \phi(t) = (\phi_1(t), \ldots, \phi_N(t))$ be an absolutely continuous mapping from $[a, b]$ to $\mathbb{R}_+^N$. When $\phi_i(t) > 0$ for all $t \in [a, b]$ and for all $i = 1, \ldots, N$, the proof of (35) is classical: for any $a \le t' < t'' \le b$ and for $\delta > 0$ small enough, the first time when the process $(X(t))$ hits the set $\bigcup_{1 \le i \le N} \{x : x_i = 0\}$ is greater than $n(t'' - t')$ whenever

$$\sup_{t \in [0, t'' - t']} |\phi(t + t') - X(nt)/n| < \delta.$$

Because of assumption $(A_3)$, this implies that

$$W_{[t', t'']}(\phi) \le -(t'' - t') L_{\{1, \ldots, N\}} \left( \frac{\phi(t'') - \phi(t')}{t'' - t'} \right)$$

and using the Markov property it follows that

$$(36) \quad W_{[a,b]}(\phi) \le \sum_{i=0}^{n-1} W_{[t_i, t_{i+1}]}(\phi) \le - \sum_{i=0}^{n-1} (t_{i+1} - t_i) L_{\{1, \ldots, N\}} \left( \frac{\phi(t_{i+1}) - \phi(t_i)}{t_{i+1} - t_i} \right)$$

for any sequence $a = t_0 < t_1 < \cdots < t_n = b$. For a piecewise linear continuous function $\phi_n : [a, b] \to \mathbb{R}_+^N$ with

$$\phi_n(t) = \phi(t_i) + (t - t_i) \frac{\phi(t_{i+1}) - \phi(t_i)}{t_{i+1} - t_i} \qquad \text{for } t \in [t_i, t_{i+1}],$$

the right-hand side of the last inequality equals

$$- \int_a^b L_{\{1,\dots,N\}}(\dot{\phi}_n(t)) \, dt.$$

When $n \to \infty$ and $\sup_i |t_{i+1} - t_i| \to 0$, $\dot{\phi}_n(t) \to \dot{\phi}(t)$ for almost all $t \in [a, b]$. By the Fatou lemma, this implies that

$$\liminf_{n \to \infty} \int_a^b L_{\{1,\dots,N\}}(\dot{\phi}_n(t)) \, dt \geq \int_a^b L_{\{1,\dots,N\}}(\dot{\phi}(t)) \, dt$$

because the convex function $L_{\{1,\dots,N\}}(\cdot)$ is lower semicontinuous. Letting therefore $n \to \infty$ and $\sup_i |t_{i+1} - t_i| \to 0$ in (36), (35) follows.

Let us prove now (35) for $\phi \in \Phi_\Lambda$ by induction with respect to $\Lambda \subset \{1, \dots, N\}$. Recall that $\Phi_\Lambda$ denotes the set of all absolutely continuous functions $\phi = (\phi_1, \dots, x_N) : [a, b] \to \mathbb{R}_+^N$ with arbitrary $a < b$, such that $\phi_i(t) > 0$ for all $i \in \Lambda$ and for all $t \in [a, b]$.

For $\phi \in \Phi_\Lambda$ with $\Lambda = \{1, \dots, N\}$, this inequality is already proved. Suppose that (35) holds for all $\phi \in \Phi_{\Lambda'}$ with $\Lambda' \subset \{1, \dots, N\}$ such that $\Lambda \subset \Lambda' \neq \Lambda$ and let us consider $\phi \in \Phi_\Lambda$, $\phi = (\phi_1, \dots, \phi_N) : [a, b] \to \mathbb{R}_+^N$. Because of our assumption, for $[t', t''] \subset [a, b]$, the inequality

$$(37) \qquad W_{[t', t'']}(\phi) \leq - \int_{t'}^{t''} L(\phi(t), \dot{\phi}(t)) \, dt$$

is already verified if there exists $i \in \{1, \dots, N\} \setminus \Lambda$ such that $\phi_i(s) > 0$ for all $s \in [t', t'']$.

Consider first the case when $\phi(t) \notin B_\Lambda$ for all $t \in [a, b]$. Then the same arguments as in the proof of Proposition 4.1 show that there is a partition $a = a_0 < a_1 < \cdots < a_n = b$ and there are $i_1, \dots, i_n \in \{1, \dots, N\} \setminus \Lambda$ such that $\phi_{i_l}(s) > 0$ for all $l = 1, \dots, n$ and for all $s \in [a_{l-1}, a_l]$. Because of our assumption, we have therefore

$$W_{[a_{l-1}, a_l]}(\phi) \leq - \int_{a_{l-1}}^{a_l} L(\phi(t), \dot{\phi}(t)) \, dt$$

for every $l = 1, \dots, n$, and hence, using the Markov property, (35) follows.

Consider now an arbitrary function $\phi \in \Phi_\Lambda$, $\phi = (\phi_1, \dots, \phi_N) : [a, b] \to \mathbb{R}_+^N$. Remark that for such a function $\phi$, (37) is already verified if $\phi(t) \notin B_\Lambda$ for all $t \in [t', t'']$. The function $\phi$ being continuous, the set

$$\Delta = \{t \in (a, b) : \phi(t) \notin B_\Lambda\} = \bigcup_{i \in \Lambda} \{t \in (a, b) : \phi_i(t) > 0\}$$

is open and consequently, it is a union of a countable collection of open disjoint intervals $(t_k, t'_k)$, $k \in \mathbb{N}$. For any $k \in \mathbb{N}$ and for any $\sigma < (t'_k - t_k)/2$ the inequality

$$W_{[t_k+\sigma, t'_k-\sigma]}(\phi) \leq - \int_{t_k+\sigma}^{t'_k-\sigma} L(\phi(t), \dot\phi(t)) \, dt$$

is therefore verified and hence, letting $\sigma \to 0$ we obtain

$$(38) \qquad W_{[t_k, t'_k]}(\phi) \leq \lim_{\sigma \to 0} W_{[t_k+\sigma, t'_k-\sigma]}(\phi) \leq - \int_{t_k}^{t'_k} L(\phi(t), \dot\phi(t)) \, dt.$$

For given $n > 0$, let us choose $k_n > 0$ large enough so that

$$(39) \qquad \sum_{k \geq k_n} t'_k - t_k < 1/n$$

and so that $k_n > k$ if $t_k = a$ or $t'_k = b$. The set $(a, b) \setminus \bigcup_{k=1}^{k_n}[t_k, t'_k]$ is a union of a finite number of disjoint intervals $(s_i, s'_i)$, $i = 1, \ldots, m$. For every $i = 1, \ldots, m$, a partition $s_i = s_{i0} < s_{i1} < \cdots < s_{ik_i} = s'_i$ is defined by induction. If $s_{ij}$ is already defined:

(a) we let $s_{ij+1} = s'_i$ and $k_i = j + 1$, if $s'_i < s_{ij} + 2/n$;

(b) otherwise, because of (38), there is $s_{ij} + 1/n < s < s_{ij} + 2/n$ such that $\phi(s) \in B_\Lambda$ and we let $s_{ij+1} = s$.

Remark that for any $i = 1, \ldots, m$ and for any $j = 1, \ldots, k_i$, by construction $\phi(s_{ij-1})$, $\phi(s_{ij}) \in B_\Lambda$. Moreover, for $\delta > 0$ small enough, the first time when the Markov process $(X(t))$ hits the set $\bigcup_{i \in \Lambda}\{x : x_i = 0\}$ is greater than $n(s_{ij} - s_{ij-1})$ whenever the inequality

$$\sup_{t \in [0, s_{ij} - s_{ij-1}]} |\phi(t + s_{ij-1}) - X(nt)/n| < \delta$$

holds and consequently, using assumption $(A_3)$ we get

$$(40) \qquad W_{[s_{ij-1}, s_{ij}]}(\phi) \leq -(s_{ij} - s_{ij-1}) L_\Lambda \left( \frac{\phi(s_{ij}) - \phi(s_{ij-1})}{s_{ij} - s_{ij-1}} \right).$$

Define a function $\psi_n : [a, b] \to \mathbb{R}_+^N$ by setting $\psi_n(t) = \phi(t)$ for $t \in \bigcup_{k=1}^{k_n}(t_k, t'_k)$, and by setting

$$\psi_n(t) = \phi(s_{ij-1}) + (t - s_{ij-1}) \frac{\phi(s_{ij}) - \phi(s_{ij-1})}{s_{ij} - s_{ij-1}}$$
$$\text{if } t \in [s_{ij-1}, s_{ij}], \ i = 1, \ldots, k_i,$$

for $t \in [s_i, s'_i]$, $i = 1, \ldots, m$. Then using the Markov property and (38), (40) we get

$$(41) \qquad W_{[a, b]}(\phi) \leq - \int_a^b L(\phi_n(t), \dot\phi_n(t)) \, dt.$$

Notice that by construction,

$$\int_{\Delta} L(\phi_n(t), \dot{\phi}_n(t)) \, dt \geq \sum_{k=1}^{k_n} \int_{t_k}^{t'_k} L(\phi(t), \dot{\phi}(t)) \, dt \to \int_{\Delta} L(\phi(t), \dot{\phi}(t)) \, dt$$

as $n \to \infty$. Moreover,

$$\int_{[a,b] \setminus \Delta} L(\phi_n(t), \dot{\phi}_n(t)) \, dt = \int_{[a,b] \setminus \Delta} L_{\Lambda}(\dot{\phi}_n(t)) \, dt$$

because by construction $\phi_n(t) \in B_{\Lambda}$ for any $t \in [a, b] \setminus \Delta$. The function $\phi$ being absolutely continuous, we have $\dot{\phi}_n(t) \to \dot{\phi}(t)$ for almost all $t \in [a, b] \setminus \Delta$. The function $L_{\Lambda}(\cdot)$ is lower semicontinuous because it is convex and hence by the Fatou lemma,

$$\liminf_{n \to \infty} \int_{[a,b] \setminus \Delta} L_{\Lambda}(\dot{\phi}_n(t)) \, dt \geq \int_{[a,b] \setminus \Delta} L_{\Lambda}(\dot{\phi}(t)) \, dt = \int_{[a,b] \setminus \Delta} L(\phi(t), \dot{\phi}(t)) \, dt.$$

Letting therefore $n \to \infty$ on the right-hand side of (41), (35) follows. $\square$

PROOF OF THEOREM 1. Let a function $\phi : [0, T] \to \mathbb{R}_+^N$ be absolutely continuous. Then using Propositions 4.1 and 4.2 it follows that

$$W_{[0,T]}(\phi) \leq -\int_0^T L(\phi(t), \dot{\phi}(t)) \, dt \leq -I_{[0,T]}(\phi).$$

Moreover, under the hypotheses $(A_1)$, the rate function $I_{[0,T]}(\cdot)$ satisfies the lower large deviation bound (1) (see the statement (c) in Theorem 4.3 of [6]) and consequently,

$$W_{[0,T]}(\phi) \geq w_{[0,T]}(\phi) \geq -I_{[0,T]}(\phi).$$

This implies that for any absolutely continuous function $\phi : [0, T] \to \mathbb{R}_+^N$,

$$w_{[0,T]}(\phi) = W_{[0,T]}(\phi) = -\int_0^T L(\phi(t), \dot{\phi}(t)) \, dt = -I_{[0,T]}(\phi)$$

and hence, Theorem 1 is proved. $\square$

**5. Proof of Theorem 2.** We begin the proof of this theorem with the following lemma.

LEMMA 5.1. *Under the hypotheses of Theorem 2, the inequality*

$$(42) \qquad\qquad I_{[0,T]}(\phi) \geq \tilde{I}_{[0,T]}(\phi)$$

*holds for every function $\phi \in D([0, T], \mathbb{R}_+^N)$.*

PROOF. Indeed, suppose that the hypotheses $(A_1)$–$(A_3)$ are satisfied and let there exist a good rate function $\tilde{I}_{[0,T]}(\cdot)$ satisfying the upper large deviation bound

$$(43) \quad \lim_{\delta \to 0} \limsup_{n \to \infty} \sup_{y \in \mathcal{E}_n : |y-x| < \delta} \frac{1}{n} \log \mathbb{P}_y(Z_n(\cdot) \in F) \leq - \inf_{\phi \in F : \phi(0)=x} \tilde{I}_{[0,T]}(\phi)$$

for any $x \in \mathbb{R}_+^N$ and for any closed set $F \subset D([0,T],\mathbb{R}_+^N)$. Then the lower large deviation bound (1) is satisfied with the rate function $I_{[0,T]}(\cdot)$ and the inequality $w_{[0,T]}(\phi) \geq -I_{[0,T]}(\phi)$ holds for all $\phi \in D([0,T],\mathbb{R}_+^N)$. Moreover, using the upper large deviation bound (43), it follows that $W_{[0,T]}(\phi) \leq -\tilde{I}_{[0,T]}(\phi)$ and consequently (42) holds.  □

The next lemma proves that $I_{[0,T]}$ is a good rate function on $D([0,T],\mathbb{R}_+^N)$ when the hypotheses of Theorem 2 are satisfied.

LEMMA 5.2. *Under the hypotheses of Theorem 2, for any compact set $V \subset \mathbb{R}_+^N$ and for any $c > 0$ the set of all functions $\phi \in D([0,T],\mathbb{R}_+^N)$ with $\phi(0) \in V$ and $I_{[0,T]}(\phi) \leq c$ is compact in $D([0,T],\mathbb{R}_+^N)$ and every function $\phi$ with $I_{[0,T]}(\phi) < \infty$ is absolutely continuous.*

PROOF. Indeed, for any compact set $V \subset \mathbb{R}_+^N$ and for any $c > 0$, the set

$$\{\phi : \phi(0) \in V, I_{[0,T]}(\phi) \leq c\}$$

is closed in $D([0,T],\mathbb{R}_+^N)$ because the rate function $I_{[0,T]}$ is lower semicontinuous. Moreover, using (42) it follows that

$$\{\phi : \phi(0) \in V, I_{[0,T]}(\phi) \leq c\} \subset \{\phi : \phi(0) \in V, \tilde{I}_{[0,T]}(\phi) \leq c\}.$$

The rate function $\tilde{I}_{[0,T]}$ being good, the set $\{\phi : \phi(0) \in V, \tilde{I}_{[0,T]}(\phi) \leq c\}$ is compact in $D([0,T],\mathbb{R}_+^N)$ and every function $\phi : [0,T] \to \mathbb{R}_+^N$ with $\tilde{I}_{[0,T]}(\phi) < \infty$ is absolutely continuous. The set $\{\phi : \phi(0) \in V, I_{[0,T]}(\phi) \leq c\}$ is therefore also compact and every function $\phi : [0,T] \to \mathbb{R}_+^N$ with $I_{[0,T]}(\phi) < \infty$ is absolutely continuous.  □

The next lemma shows that the rate function $I_{[0,T]}$ satisfies the upper large deviation bound (2).

LEMMA 5.3. *Under the hypotheses of Theorem 2,*

$$(44) \quad \lim_{\delta \to 0} \limsup_{n \to \infty} \sup_{y \in \mathcal{E}_n : |y-x| < \delta} \frac{1}{n} \log \mathbb{P}_y(Z_n(\cdot) \in F) \leq - \inf_{\phi \in F : \phi(0)=x} I_{[0,T]}(\phi)$$

*for any $x \in \mathbb{R}_+^N$ and for any closed set $F \subset D([0,T],\mathbb{R}_+^N)$.*

PROOF.   Consider a closed set $F \subset D([0,T], \mathbb{R}_+^N)$, $x \in \mathbb{R}_+^N$ and let

$$c = \inf_{\phi \in F : \phi(0) = x} I_{[0,T]}(\phi).$$

The rate function $\tilde{I}_{[0,T]}(\cdot)$ being good, every function $\phi \in F$ satisfying the inequality $\tilde{I}_{[0,T]}(\phi) \leq c$ is absolutely continuous and hence, by Theorem 1,

$$W_{[0,T]}(\phi) \leq -I_{[0,T]}(\phi).$$

The last inequality implies that for any $\epsilon > 0$ there exists $\delta_\phi > 0$ such that

$$\limsup_{n \to \infty} \sup_{z \in \mathcal{E}_n : |z - \phi(0)| < \delta_\phi} \frac{1}{n} \log \mathbb{P}_z(\|\phi - Z_n\| < \delta_\phi) \leq -I_{[0,T]}(\phi) + \epsilon$$

for all $0 < \delta < \delta_\phi$. The set $K = \{\phi \in F : \phi(0) = x, \tilde{I}_{[0,T]}(\phi) \leq c\}$ being compact, there exists a finite collection of functions $\phi_1, \ldots, \phi_n \in K$ such that

$$K \subset \mathcal{O} = \bigcup_{i=1}^{n} \{\phi : \|\phi_i - \phi\| < \delta_{\phi_i}\}$$

and consequently,

(45)   $$\lim_{\delta \to 0} \limsup_{n \to \infty} \sup_{y \in \mathcal{E}_n : |y - x| < \delta} \frac{1}{n} \log \mathbb{P}_y(Z_n(\cdot) \in \mathcal{O}) \leq -\inf_{1 \leq i \leq n} I_{[0,T]}(\phi_i) + \epsilon.$$

Remark that for any $\phi \in F \setminus \mathcal{O}$ with $\phi(0) = x$, the following inequality holds:

$$I_{[0,T]}(\phi) \geq \tilde{I}_{[0,T]}(\phi) \geq c = \inf_{\phi \in F : \phi(0) = x} I_{[0,T]}(\phi).$$

The set $F \setminus \mathcal{O}$ being closed, using the upper large deviation bound (43), this implies that

$$\lim_{\delta \to 0} \limsup_{n \to \infty} \sup_{y \in \mathcal{E}_n : |y - x| < \delta} \frac{1}{n} \log \mathbb{P}_y(Z_n(\cdot) \in F \setminus \mathcal{O}) \leq -\inf_{\phi \in F : \phi(0) = x} I_{[0,T]}(\phi).$$

This inequality and (45) show that

$$\lim_{\delta \to 0} \limsup_{n \to \infty} \sup_{y \in \mathcal{E}_n : |y - x| < \delta} \frac{1}{n} \log \mathbb{P}_y(Z_n(\cdot) \in F) \leq -\inf_{\phi \in F : \phi(0) = x} I_{[0,T]}(\phi) + \epsilon.$$

Letting in the last inequality $\varepsilon \to 0$, (44) follows.   □

The last lemma completes the proof of Theorem 2 (the lower large deviation bound (1) is satisfied because of assumption $(A_1)$; this is a consequence of the statement (c) of Theorem 4.3 of [6]).

## APPENDIX

In this section we prove inequality (11). This is a subject of the following lemma.

LEMMA A.1.  *Under the hypotheses* $(B_0)$–$(B_2)$, *for any* $\Lambda \subset \{1, \ldots, N\}$, $t > 0$ *and* $v \in \mathbb{R}^\Lambda$, *the local Markov-additive process* $(A_\Lambda(t), Y_\Lambda(t))$ *satisfies the inequality*

$$(46) \quad \lim_{\delta \to 0} \liminf_{n \to \infty} \frac{1}{n} \log \mathbb{P}_{(0,0)} \left( \sup_{s \in [0, nt]} |A_\Lambda(s) - sv| + |Y_\Lambda(s)| < \delta n \right) \geq -t \lambda_\Lambda^*(v).$$

PROOF.   Given $K \subset \mathbb{Z}_+^{\Lambda^c}$, let $T_K$ be the first time when the process $(Y_\Lambda(t))$ exists from the set $K$ and let $\mathcal{K}_\Lambda$ be the collection of all the finite subsets $K$ of $\mathbb{Z}_+^{\Lambda^c}$ for which the restriction of the Markov chain $(Y_\Lambda(t))$ on $K$ is irreducible. For $K \subset \mathcal{K}_\Lambda$, the matrices $\mathcal{P}_{\Lambda,K}^t(\alpha) = (\mathcal{P}_{\Lambda,K}^t(\alpha; y, y')$; $y, y' \in K)$ with

$$\mathcal{P}_{\Lambda,K}^t(\alpha; y, y') = \mathbb{E}_{(0,y)}(e^{\langle \alpha, A_\Lambda(t) \rangle} \mathbb{1}_{\{Y_{\Lambda^c}(t) = y' \text{ and } T_K > t\}})$$

are irreducible. Moreover, $\mathcal{P}_{\Lambda,K}^t(\alpha) = \exp(t \cdot Q_{\Lambda,K}(\alpha))$ where the matrix $Q_{\Lambda,K}(\alpha) = (Q_{\Lambda,K}(\alpha; y, y')$; $y, y' \in K)$ is defined by

$$Q_{\Lambda,K}(\alpha; y, y') = \sum_{x \in \mathbb{Z}^\Lambda} q_\Lambda((0, y), (x, y')) e^{\langle \alpha, x \rangle}.$$

Using the Perron–Frobenius theorem this implies that:

(a) The matrix $Q_{\Lambda,K}(\alpha)$ has a unique maximal real eigenvalue $\lambda_{\Lambda,K}(\alpha)$ and a strictly positive unique to constant multiples right eigenvector $f_{\Lambda,K}^\alpha = (f_{\Lambda,K}^\alpha(y)$; $y \in K)$ associated with $\lambda_{\Lambda,K}(\alpha)$.

(b) For every $t > 0$, $r_{\Lambda,K}^t(\alpha) = \exp\{t \lambda_{\Lambda,K}(\alpha)\}$ is the unique Perron–Frobenius eigenvalue of the matrix $\mathcal{P}_{\Lambda,K}^t(\alpha)$ and $f_{\Lambda,K}^\alpha = (f_{\Lambda,K}^\alpha(y)$; $y \in K)$ is its unique to constant multiples right eigenvector associated with $r_{\Lambda,K}^t(\alpha)$.

(c) The collection of the functions $\lambda_{\Lambda,K}$, $K \in \mathcal{K}_\Lambda$, is increasing with respect to $K$ and for all $y, y' \in K$,

$$\lambda_{\Lambda,K}(\alpha) = \limsup_{t \to \infty} \frac{1}{t} \log \mathcal{P}_{\Lambda,K}^t(\alpha; y, y') \leq \lambda_\Lambda(\alpha).$$

Moreover, an argument similar to one used to prove Lemma 1 of [11] shows that the functions $\lambda_{\Lambda,K}$, $K \in \mathcal{K}_\Lambda$, are convex and infinitely differentiable on $\mathbb{R}^\Lambda$, and using the arguments of the proof of Proposition 2 in [11] we obtain that

$$\lambda_\Lambda(\alpha) = \sup_{K \in \mathcal{K}_\Lambda} \lambda_{\Lambda,K}(\alpha).$$

Let $\lambda^*_{\Lambda,K}$ be the convex conjugate of the function $\lambda_{\Lambda,K}$. Then the collection of the functions $\lambda^*_{\Lambda,K}$, $K \in \mathcal{K}_\Lambda$, is decreasing with respect to $K$ and using Theorem 16.5 of [12] it follows that the convex conjugate $\lambda^*_\Lambda$ of the function $\lambda_\Lambda$ is the closure of the function

$$\inf_{K \in \mathcal{K}_\Lambda} \lambda^*_{\Lambda,K}.$$

The mapping

$$v \to \lim_{\delta \to 0} \liminf_{n \to \infty} \frac{1}{n} \log \mathbb{P}_{(0,0)} \left( \sup_{s \in [0,nt]} |A_\Lambda(s) - sv| + |Y_\Lambda(s)| < \delta n \right)$$

being upper semicontinuous on $\mathbb{R}^\Lambda$, we conclude that to prove (46) it is sufficient to show that

$$\lim_{\delta \to 0} \liminf_{n \to \infty} \frac{1}{n} \log \mathbb{P}_{(0,0)} \left( \sup_{s \in [0,nt]} |A_\Lambda(s) - sv| + |Y_\Lambda(s)| < \delta n \right) \geq -t \inf_{K \in \mathcal{K}_\Lambda} \lambda^*_{\Lambda,K}(v).$$

To prove the last inequality, it is sufficient to show that for any finite set $K \in \mathcal{K}_\Lambda$,

$$(47) \quad \lim_{\delta \to 0} \liminf_{n \to \infty} \frac{1}{n} \log \mathbb{P}_{(0,0)} \left( \sup_{s \in [0,nt]} |A_\Lambda(s) - sv| < \delta n, T_K > nt \right) \geq -t\lambda^*_{\Lambda,K}(v).$$

For this we use the martingale method and the classical method of change of measure.

Let $\{\mathcal{F}^\Lambda_t\}_{t \geq 0}$ be the natural filtration of the Markov process $(A_\Lambda(t), Y_\Lambda(t))$. For $(x,y) \in \mathbb{Z}^\Lambda \times K$ and $t \geq 0$ we define a new measure $\mathbb{P}^{(\alpha)}_{(x,y)}$ on

$$\mathcal{F}^\Lambda_t \cap \{T_K > t\} = \{E \cap \{T_K > t\} : E \in \mathcal{F}^\Lambda_t\}$$

by setting

$$\mathbb{P}^{(\alpha)}_{(x,y)}(B) = \mathbb{E}_{(x,y)}(\mathbb{1}_B \exp\{\langle \alpha, A_\Lambda(t) - x \rangle - t\lambda_{\Lambda,K}(\alpha)\} f^\alpha_{\Lambda,K}(Y_\Lambda(t))/f^\alpha_{\Lambda,K}(y)).$$

Then clearly $\mathbb{P}^{(\alpha)}_{(x,y)}(T_K > t) = 1$ for all $(x,y) \in \mathbb{Z}^\Lambda \times K$ and for all $t > 0$. This implies that for all $s > t > 0$ and for all $E \in \mathcal{F}^\Lambda_t$,

$$\mathbb{P}^{(\alpha)}_{(x,y)}(E \cap \{T_K > s\}) = \mathbb{P}^{(\alpha)}_{(x,y)}(E \cap \{T_K > t\})$$

and hence, letting

$$\mathbb{P}^{(\alpha)}_{(x_0,y_0)}(E \cap \{T_K = \infty\}) = \lim_{s \to \infty} \mathbb{P}^{(\alpha)}_{(x_0,y_0)}(E \cap \{T_K > s\})$$

for $E \in \bigcup_{t \geq 0} \mathcal{F}^\Lambda_t$ we obtain a new probability measure on $\bigcup_{t \geq 0} \mathcal{F}^\Lambda_t \cap \{T_K = \infty\}$. This is a distribution of a new Markov process on $\mathbb{Z}^\Lambda \times K$ with initial state $(x_0, y_0)$ and transition intensities

$$q^{(\alpha)}_{\Lambda,K}((x,y),(x',y')) = q_\Lambda((x,y),(x',y'))e^{\langle \alpha, x'-x \rangle} f^\alpha_{\Lambda,K}(y')/f^\alpha_{\Lambda,K}(y)$$

$x, x' \in \mathbb{Z}^\Lambda$, $y, y' \in K$.

Let $\mathbb{E}^{(\alpha)}_{(x,y)}$ denote the expectation with respect to the new probability measure $\mathbb{P}^{(\alpha)}_{(x,y)}$ and let $E_{n,\delta} = \{\sup_{s \in [0,nt]} |A_\Lambda(s) - sv| < \delta n, T_K > nt\}$. Without any restriction of generality we will suppose that $0 \in K$ and that $f^\alpha_{\Lambda,K}(0) = 1$. Then using the standard arguments of the change of measure it follows that

$$\log \mathbb{P}_{(0,0)}(E_{n,\delta})$$
$$= \log \mathbb{E}^{(\alpha)}_{(0,0)}(\mathbb{1}_{E_{n,\delta}} \exp\{-\langle \alpha, A_\Lambda(nt) \rangle + nt\lambda_{\Lambda,K}(\alpha)\} f^\alpha_{\Lambda,K}(Y_\Lambda(nt)))$$
$$\geq \log \mathbb{P}^{(\alpha)}_{(0,0)}(E_{n,\delta}) + nt(\lambda_{\Lambda,K}(\alpha) + \langle \alpha, v \rangle) - \delta n - \max_{y \in K} \log f^\alpha_{\Lambda,K}(y).$$

Suppose first that $v \in \mathrm{ri}(\mathrm{dom}\,\lambda^*_K)$. Then, the function $\lambda_K(\cdot)$ being convex and differentiable on $\mathbb{R}^\Lambda$, there is $\alpha_v \in \mathbb{R}^\Lambda$ such that

$$(48) \qquad \lambda^*_{\Lambda,K}(v) = \langle \alpha_v, v \rangle - \lambda_{\Lambda,K}(\alpha_v)$$

(see Corollary 26.4.1 of [12]). Using the last inequality with $\alpha = \alpha_v$ we obtain

$$\lim_{\delta \to 0} \liminf_{n \to \infty} \frac{1}{n} \log \mathbb{P}_{(0,0)}(E_{n,\delta}) \geq \lim_{\delta \to 0} \liminf_{n \to \infty} \frac{1}{n} \log \mathbb{P}^{(\alpha_v)}_{(0,0)}(E_{n,\delta}) + \lambda^*_{\Lambda,K}(v)$$

and hence, to get (47) it is sufficient to show that for any $\delta > 0$,

$$\mathbb{P}^{(\alpha_v)}_{(0,0)}(E_{n,\delta}) \to 1 \qquad \text{as } n \to \infty$$

or equivalently that

$$(49) \qquad \mathbb{P}^{(\alpha_v)}_{(0,0)}\left(\sup_{s \in [0,nt]} |A_\Lambda(s) - sv| \geq \delta n\right) \to 0 \qquad \text{as } n \to \infty.$$

For this we use a martingale technique. Straightforward calculations show that for any $\alpha \in \mathbb{R}^\Lambda$,

$$M(\alpha, t) = \mathbb{1}_{\{T_K > t\}} \exp\{\langle \alpha - \alpha_v, A_\Lambda(t) \rangle - (\lambda_{\Lambda,K}(\alpha) - \lambda_{\Lambda,K}(\alpha_v))t\}$$
$$\times f^\alpha_{\Lambda,K}(Y_\Lambda(t))/f^{\alpha_v}_{\Lambda,K}(Y_\Lambda(t))$$

is a martingale relative to the new probability measure $\mathbb{P}^{(\alpha_v)}_{(0,0)}$ with

$$\mathbb{E}^{(\alpha_v)}_{(0,0)}(M(\alpha, t)) \equiv 1.$$

Moreover, because of (48)

$$M(\alpha, t) = \mathbb{1}_{\{T_K > t\}} \exp\{\langle \alpha - \alpha_v, A_\Lambda(t) - vt \rangle - (\lambda_{\Lambda,K}(\alpha) - \langle \alpha, v \rangle + \lambda^*_{\Lambda,K}(v))t\}$$
$$\times f^\alpha_{\Lambda,K}(Y_\Lambda(t))/f^{\alpha_v}_{\Lambda,K}(Y_\Lambda(t))$$

and hence, using Fenchel's inequality $\lambda_{\Lambda,K}(\alpha) - \langle \alpha, v \rangle + \lambda_{\Lambda,K}^*(v) \geq 0$, it follows that

$$\Xi(t) = M(\alpha, t) \exp\{(\lambda_{\Lambda,K}(\alpha) - \langle \alpha, v \rangle + \lambda_{\Lambda,K}^*(v))t\}$$
$$= \exp\{\langle \alpha - \alpha_v, X(t) - vt \rangle\} f_{\Lambda,K}^\alpha(Y_\Lambda(t))/f_{\Lambda,K}^{\alpha_v}(Y_\Lambda(t))$$

is a submartingale relative to the new probability measure $\mathbb{P}_{(0,0)}^{(\alpha_v)}$ with

$$\mathbb{E}_{(0,0)}^{(\alpha_v)}(\Xi(t)) = \exp\{t(\lambda_{\Lambda,K}(\alpha) - \langle \alpha, v \rangle + \lambda_{\Lambda,K}^*(v))\}.$$

Letting $c(\alpha) = \min_{y \in K} f_{\Lambda,K}^\alpha(y)/f_{\Lambda,K}^{\alpha_v}(y)$ and using submartingale inequality, it follows that for any $\gamma > 0$, and for any $\alpha \in \mathbb{R}^\Lambda$ with $|\alpha - \alpha_v| \leq 1$,

$$\mathbb{P}_{(0,0)}^{(\alpha_v)}\left( \sup_{s \in [0,nt]} \langle \alpha - \alpha_v, A_\Lambda(s) - vs \rangle \geq \gamma \right)$$

(50)
$$\leq \mathbb{P}_{(0,0)}^{(\alpha_v)}\left( \sup_{s \in [0,nt]} \Xi(s) \geq c(\alpha)e^\gamma \right)$$

$$\leq c^{-1}(\alpha) \exp(-\gamma + nt(\lambda_{\Lambda,K}(\alpha) - \langle \alpha, v \rangle + \lambda_{\Lambda,K}^*(v))).$$

Moreover, let

$$C = \max_{\alpha:\, |\alpha - \alpha_v| \leq 1} \|\partial_\alpha^2 \lambda_{\Lambda,K}(\alpha)\|,$$

where $\partial_\alpha^2 \lambda_{\Lambda,K}(\alpha)$ denotes Hessian matrix of $\lambda_{\Lambda,K}(\cdot)$. Then for any $\alpha \in \mathbb{R}^\Lambda$ with $|\alpha - \alpha_v| \leq 1$, the inequality $\lambda_{\Lambda,K}(\alpha) - \langle \alpha, v \rangle + \lambda_{\Lambda,K}^*(v) \leq |\alpha - \alpha_v|^2 C$ holds and using (50) we get

$$\mathbb{P}_{(0,0)}^{(\alpha_v)}\left( \sup_{s \in [0,nt]} \langle \alpha - \alpha_v, A_\Lambda(s) - vs \rangle \geq \gamma \right)$$

$$\leq c^{-1}(\alpha) \exp(-\gamma + nt|\alpha - \alpha_v|^2 C).$$

Let $\epsilon \in \mathbb{R}\Lambda$ be a unit vector. Letting in the above inequality $\alpha = \alpha_v + \theta\epsilon$ with $0 < \theta < 1$ and $\gamma = (Ct+1)\theta^2 n$, we obtain

$$\mathbb{P}_{(0,0)}^{(\alpha_v)}\left( \sup_{s \in [0,nt]} \langle \epsilon, A_\Lambda(s) - vs \rangle \geq (Ct+1)\theta n \right) \leq c^{-1}(\alpha) \exp(-\theta^2 n).$$

Finally, the unit vector $\epsilon$ being arbitrary, the last inequality proves that

$$\mathbb{P}\left( \sup_{s \in [0,nt]} |A_\Lambda(s) - vs| \geq 2N\theta(C\tau+1)n \right)$$

$$\leq 2N \max_\epsilon c^{-1}(\alpha_v + \theta\epsilon) \exp\{-\theta^2 n\}$$

and hence, letting $\delta = 2N\theta(C\tau+1)$, (49) follows.

For $v \in \mathrm{ri}(\mathrm{dom}\,\lambda_K^*)$, (47) is therefore verified. The mapping

$$v \to \lim_{\delta \to 0} \liminf_{n \to \infty} \frac{1}{n} \log \mathbb{P}_{(0,0)} \left( \sup_{s \in [0,nt]} |A_\Lambda(s) - sv| < \delta n, T_K > nt \right)$$

being upper semicontinuous on $\mathbb{R}^\Lambda$, this implies that (47) holds for every $v \in \mathbb{R}^\Lambda$ and hence, Lemma A.1 is proved. $\quad\square$

DÉPARTEMENT DE MATHÉMATIQUES
UNIVERSITÉ DE CERGY-PONTOISE
2 AVENUE ADOLPHE CHAUVIN
95302 CERGY-PONTOISE CEDEX
FRANCE
E-MAIL: Irina.Ignatiouk@math.u-cergy.fr